\documentclass[11pt]{amsart}

\usepackage[margin=1.12in]{geometry}
\usepackage{amsmath,amssymb,mathtools}
\usepackage{enumitem}
\usepackage{xcolor}
\usepackage{microtype}
\usepackage[colorlinks=true,linkcolor=blue!55!black,
citecolor=blue!55!black,urlcolor=blue!55!black]{hyperref}
\usepackage{microtype}
\usepackage{needspace}
\usepackage{xcolor}
\newtheorem{theorem}{Theorem}[section]
\newtheorem{proposition}[theorem]{Proposition}
\newtheorem{corollary}[theorem]{Corollary}
\newtheorem{lemma}[theorem]{Lemma}
\theoremstyle{remark}

\numberwithin{equation}{section}
\DeclareMathOperator{\PSH}{PSH}
\allowdisplaybreaks
\renewenvironment{abstract}
  {\small\begin{quote}\noindent\textsc{\abstractname}.\enspace\ignorespaces}
  {\end{quote}}
\title[The Dirichlet problem for complex Monge--Amp\`ere equations]
{The Dirichlet problem for complex Monge--Amp\`ere equations
with H\"older continuous right-hand sides}
\author{Liding Huang}
\address{School of Mathematical Sciences, Xiamen University,
Xiamen 361005, P. R. China}
\email{huangliding@xmu.edu.cn}
\author{Xi Zhang}
\address{School of Mathematics and Statistics, Nanjing University of Science and Technology,
Nanjing 210094, P. R. China}
\email{mathzx@njust.edu.cn}
\date{}
\begin{document}

\maketitle
\begingroup
\renewcommand{\thefootnote}{}
\endgroup

\begin{abstract}
We study the Dirichlet problem for the complex Monge--Amp\`ere
equation on bounded strongly pseudoconvex domains with $C^3$ boundary
and $C^3$ boundary values. For every $n\ge2$ and $0<\alpha<1$,
each continuous Bedford--Taylor solution with a strictly positive
$C^\alpha$ right-hand side satisfies a global $C^{2,\alpha}$ estimate.
There is a unique such solution, which is strictly plurisubharmonic
up to the boundary.
Combining new boundary estimates with weighted integral estimates
for the trace of the complex Hessian along an auxiliary parabolic flow, we obtain
global $C^{1,\beta}$ bounds for every fixed $0<\beta<1$ in every
dimension. Green function estimates and
comparison arguments near the boundary give the boundary Hessian
bound. Solutions of auxiliary Dirichlet problems with locally
constant right-hand sides have uniform $W^{2,2n^2}$ bounds.
Rescaling and interpolation then reduce the global Hessian estimate
to interior regularity.
\end{abstract}


\section{Introduction}\label{sec:introduction}

Let $\Omega\Subset\mathbb C^n$ be a bounded strongly pseudoconvex
domain with $C^3$ boundary. We study the global regularity of solutions
to the following Dirichlet problem for the complex Monge--Amp\`ere equation
\begin{equation}\label{intro:dirichlet}
 \begin{cases}
  \det(u_{i\bar j})=f&\text{in }\Omega,\\
  u=\varphi&\text{on }\partial\Omega,
 \end{cases}
 \qquad u\in\PSH(\Omega)\cap C^0(\overline\Omega),
\end{equation}
where $\PSH(\Omega)$ denotes the set of real-valued plurisubharmonic
functions on $\Omega$, $\varphi\in C^3(\overline\Omega)$ is real-valued,
and $f\in C(\overline\Omega)$ is strictly positive. The first
equation is understood as the Bedford--Taylor measure identity
\[
 (i\partial\bar\partial u)^n=f(i\partial\bar\partial |z|^2)^n
 \quad\text{in }\Omega.
\]

Bedford--Taylor
\cite{BT} established existence and uniqueness of
the continuous plurisubharmonic solution on strongly pseudoconvex
domains. Caffarelli--Kohn--Nirenberg--Spruck
\cite{CKNS} subsequently established classical solvability for smooth
positive right-hand sides and smooth boundary data on smooth strongly
pseudoconvex domains.
For $f\ge0$, Krylov \cite{Krylov} established global $C^{1,1}$
regularity when the boundary and boundary values are of class
$C^{3,1}$ and $f^{1/n}\in C^{1,1}(\bar\Omega)$. Dinew--Sroka \cite{DS} later gave an analytic proof of this result.
For smooth data on a general smooth bounded domain,
Guan \cite{Guan} replaced strong pseudoconvexity by the existence
of a strictly plurisubharmonic subsolution with the prescribed boundary value.
For $f\in L^p$, $p>1$, global H\"older regularity was studied in
\cite{BKPZ,Charabati,GKZ,HZ}.

For the real Monge--Amp\`ere equation, Wang \cite{XW96}
obtained boundary regularity on strictly convex $C^3$ domains with
$C^3$ boundary values and $f\in C^{1,1}$. Trudinger--Wang \cite{TW} proved the global $C^{2,\alpha}$ estimate on
uniformly convex $C^3$ domains with $C^3$ boundary values and a positive
$C^\alpha$ right-hand side, for every $0<\alpha<1$. Savin \cite{Savin} developed localization of sections at the boundary
and pointwise $C^{2,\alpha}$ estimates for the real equation. The pointwise estimates require $C^{2,\alpha}$ boundary geometry
and data, together with quadratic separation from the supporting plane.
On uniformly convex domains, Savin also obtained estimates under
pointwise $C^3$ assumptions on the boundary geometry and data, and
global estimates for zero boundary values on $C^{2,\alpha}$ domains.

A natural question is whether the global boundary estimates of
Trudinger--Wang \cite{TW} extend to the complex Monge--Amp\`ere
equation. Chen--He \cite{CH} asked whether a strictly positive
$C^\alpha$ right-hand side and
zero boundary values on a bounded strongly pseudoconvex domain imply
$C^{2,\alpha}$ regularity of the solution. We study the Dirichlet problem with $C^3$ geometry and arbitrary $C^3$ boundary
values, for every $0<\alpha<1$ and every complex dimension $n\ge2$.

For the Dirichlet problem with zero boundary values on bounded
convex domains, B{\l}ocki \cite{Blocki} proved interior smoothness
for smooth positive $f$ with $D(f^{1/n})\in L^\infty$.
Further interior estimates use a priori regularity of the solution.
Dinew--Zhang--Zhang
\cite{DZZ} obtained interior $C^{2,\alpha}$ regularity from
$C^{1,1}$ control.
Assuming a strictly positive $C^{1,1}$ right-hand side,
B{\l}ocki--Dinew \cite{BD} proved interior Laplacian bounds under
$W^{2,p}$ control for $p>n(n-1)$. He \cite{He} derived interior
Laplacian bounds from $L^p$ control of the Laplacian for $p>n^2$,
under a Lipschitz bound for $\log f$.
Under a Laplacian bound, Wang \cite{YW} reduced the equation to a
uniformly elliptic one and obtained
$C^{2,\gamma}$ regularity for every $\gamma<\alpha$. Li--Li--Zhang \cite{LLZ} proved the
interior $C^{2,\alpha}$ estimate assuming $C^{1,\beta}$ regularity
with $\beta>1-\alpha/[n(2+\alpha)-1]$.
Hu--Sheng \cite{HS} obtained interior $C^2$ estimates for positive
Dini-continuous densities under quadratic growth bounds and a
distributional lower bound for $\Delta(f^{1/n})$.

For positive H\"older right-hand sides, the gradient estimate
of Chen--Cheng \cite{CC} applies to the
Dirichlet problem once global $C^0$ and boundary gradient bounds are
available. The barriers below provide these bounds. Cheng--Xu \cite{CX} established
perturbative interior $W^{2,p}$ estimates, including estimates near a
smooth strictly plurisubharmonic function. 
Li--Wang--Zheng \cite{LWZ}
proved the global $C^{2,\alpha}$ estimate for every $0<\alpha<1$ on
compact K\"ahler manifolds without boundary.

The following theorem establishes the global $C^{2,\alpha}$
estimate up to the boundary on strongly pseudoconvex domains.

\begin{theorem}\label{thm:main}
Let $n\geq2$, let $\Omega\Subset\mathbb C^n$ be a bounded strongly
pseudoconvex domain with $C^3$ boundary, and let
$\varphi\in C^3(\bar\Omega)$ be real-valued. Let $0<\alpha<1$ and suppose
\[
 0<\lambda\leq f\leq\Lambda,\qquad
 \|f\|_{C^\alpha(\overline\Omega)}\leq M.
\]
If $u\in\PSH(\Omega)\cap C^0(\bar\Omega)$ solves
\eqref{intro:dirichlet} in the Bedford--Taylor sense, then
$u\in C^{2,\alpha}(\bar\Omega)$ and
\begin{equation}
 \|u\|_{C^{2,\alpha}(\overline\Omega)}
 \le C,\label{eq:SE2}
\end{equation}
where $C$ depends only on $n,\alpha,\lambda,\Lambda,M$, the quantitative
$C^3$ geometric bounds for $\Omega$, the $C^3$ norm of $\varphi$,
and the continuity moduli of the third derivatives of the boundary
charts and $\varphi$. These third derivatives need only be continuous.
\end{theorem}

We also obtain classical solvability under these assumptions.

\begin{theorem}\label{thm:existence}
Let $n\ge2$, let $\Omega\Subset\mathbb C^n$ be a bounded strongly
pseudoconvex domain with $C^3$ boundary, and let
$\varphi\in C^3(\bar\Omega)$ be real-valued. Let $0<\alpha<1$. If
\[
 f\in C^\alpha(\bar\Omega),\qquad
                    \inf_{\bar\Omega}f>0,
\]
then \eqref{intro:dirichlet} has a unique continuous
plurisubharmonic solution $u$. Moreover, $u$ belongs to
$C^{2,\alpha}(\bar\Omega)$ and is strictly plurisubharmonic up to the
boundary.

\end{theorem}

\Needspace{5\baselineskip}
For smoother data, the usual elliptic bootstrap gives the following result.

\begin{corollary}\label{cor:higher-regularity}
Let $n\ge2$, $k\in\mathbb N$, $k\ge1$, and $0<\sigma<1$.
Suppose $\Omega\Subset\mathbb C^n$ has strongly pseudoconvex
$C^{k+2,\sigma}$ boundary, $\varphi\in C^{k+2,\sigma}(\bar\Omega)$
is real-valued, and $f\in C^{k,\sigma}(\bar\Omega)$ is strictly positive.
Then the solution $u$ of \eqref{intro:dirichlet} given by
Theorem~\ref{thm:existence} satisfies
\[
 \|u\|_{C^{k+2,\sigma}(\bar\Omega)}\le C,
\]
where $C$ depends only on $n,k,\sigma,\Omega$,
$\|\varphi\|_{C^{k+2,\sigma}}$, $\|f\|_{C^{k,\sigma}}$ and
$(\inf_{\bar\Omega}f)^{-1}$.
\end{corollary}

Two difficulties arise in the boundary second derivative estimate
for Theorem~\ref{thm:main}. First, we need a uniform positive lower
bound for the complex tangential Hessian under $C^3$ assumptions
on the boundary geometry and boundary values. Second, the mixed
derivatives must be controlled using only the $C^\alpha$ norm of $f$,
whereas the classical estimates \cite{CKNS,Guan2009} require bounds
for $Df$.

For the tangential lower bound, we regularize a $C^1$ coefficient in
the comparison function and adjust its quadratic normal term.
By choosing the smoothing scale and the correction together, we
make the determinant error tend to zero using the $C^3$ bounds
and the continuity moduli of third derivatives of the boundary charts
and boundary values. The corresponding construction in
\cite{Guan2009} uses a $C^4$ defining function.

For the mixed derivatives, we need gradient H\"older exponents close
to one. Starting from the gradient bound of Chen--Cheng \cite{CC}, we
modify the elliptic and parabolic approximation method of
Li--Wang--Zheng \cite{LWZ}. Their energy estimate does not directly
control the boundary flux produced by integration by parts.
Along the auxiliary parabolic flow, we introduce an additional weight
comparable to the boundary distance. A separate boundary estimate
bounds this flux by an arbitrarily small multiple of the weighted
interior dissipation, plus a constant, so it can be absorbed. The resulting weighted trace bounds give gradient
estimates at positive times. Interpolation with the $C^0$
approximation error and iteration yield global $C^{1,\beta}$ bounds
for every fixed $0<\beta<1$.

By Young's inequality for Stieltjes integrals
\cite{Young}, we control the
heat-regularized products of tangential derivatives of $f$ with
$C^\beta$ coefficients using only H\"older norms. The $C^{1,\beta}$
bound also gives an $L^1$ gradient estimate for the Dirichlet Green
function that decays as the pole approaches the boundary.
Together with comparison functions whose gradients are $O(R^\beta)$
on neighborhoods of radius $R$, this decay gives the mixed derivative
bounds when $\alpha+\beta/(2-\beta)>1$. Choosing $\beta$ close to one
covers every $0<\alpha<1$. The tangential lower bound and determinant
equation then give the full real boundary Hessian bound.

For the global Hessian bound, we combine coefficient integrability
with the boundary gradient modulus. By quantifying the comparison
argument of Li--Li--Zhang \cite{LLZ}, we obtain, for $\beta$ sufficiently
close to one,
\[
 |D^2u(x)|\le C\,\operatorname{dist}(x,\partial\Omega)^{-(n+1)(1-\beta)}.
\]
This yields uniform $W^{2,2n^2}$ bounds and hence the coefficient
integrability needed for Moser iteration.

We then solve an auxiliary Dirichlet problem whose density is
constant near a boundary point. The second directional derivatives
of its solution are subsolutions of the linearized equation there. By Moser iteration
and the comparison functions above, we obtain a uniform Hessian bound
in a thin boundary neighborhood. The thickness is chosen so that
rescaling by the boundary distance gives a uniform $C^{1,1}$ bound
for the auxiliary solution and a vanishing comparison error.
Interpolation using the known gradient H\"older seminorms gives a
uniform $C^{1,1-\alpha/(4n)}$ bound for the rescaled original solution.
The global Hessian bound follows from \cite{LLZ}.
Uniformly elliptic boundary regularity
\cite{SS} and the Schauder estimates
complete the $C^{2,\alpha}$ estimate.

The paper is organized as follows. Section~\ref{sec:preliminaries} contains the geometric constructions
and auxiliary estimates. The gradient and boundary Hessian estimates
are proved in Section~\ref{sec:gradient-boundary}, followed by the
global Hessian bound in Section~\ref{sec:global-hessian}.
Section~\ref{sec:main-proof} completes the proofs of the main results.

\Needspace{5\baselineskip}
\medskip
\noindent\textbf{Acknowledgements.} The authors gratefully acknowledge
 the assistance of AI models, particularly ChatGPT 6.1 Astra, in exploring
proof strategies and carrying out calculations. All AI-generated suggestions
have been verified and revised by the authors, who take full responsibility
for all assertions in this paper.

\section{Boundary geometry and auxiliary estimates}\label{sec:preliminaries}

In this section, we construct barriers and smooth approximations,
prove the lower bound for the complex tangential Hessian, and
establish estimates uniform under smooth approximation.

For $z_j=x_j+i y_j$, we use the conventions
\[
 \partial_i=\partial_{z_i}=\tfrac12(\partial_{x_i}-i\partial_{y_i}),
 \qquad
 \partial_{\bar i}=\partial_{\bar z_i}
     =\tfrac12(\partial_{x_i}+i\partial_{y_i}).
\]
Throughout, $D$ denotes differentiation in real coordinates.
For real-valued $u$, set
\[
 \begin{gathered}
  u_i=\partial_i u,\qquad
  u_{i\bar j}=\partial_i\partial_{\bar j}u,\\
  \partial u=(u_i),\qquad |\partial u|_0^2=\sum_i|u_i|^2.
 \end{gathered}
\]
Repeated complex indices are summed from $1$ to $n$.
For a strictly plurisubharmonic function $u$, we write
\begin{equation}
 u^{i\bar j}=\bigl[(u_{k\bar\ell})^{-1}\bigr]_{ji},
 \qquad \sum_j u^{i\bar j}u_{k\bar j}=\delta_{ik}.
 \label{eq:inverse-coefficients}
\end{equation}

For real-valued $h$, write
\[
 L_uh=u^{i\bar j}h_{i\bar j},\qquad
 |\partial h|_u^2=u^{i\bar j}h_i h_{\bar j}.
\]%
For smooth strictly plurisubharmonic $u$, we obtain the cofactor
identities by differentiating the determinant and commuting derivatives:
\begin{equation}
 C^{i\bar j}=fu^{i\bar j},\qquad f=\det(u_{k\bar\ell}),\qquad
 \sum_j\partial_{\bar j}C^{i\bar j}=0,\qquad
 \sum_i\partial_iC^{i\bar j}=0.
 \label{eq:cofactor-divergence}
\end{equation}
The real Euclidean Laplacian and the standard complex structure satisfy
\[
 \Delta u=4\sum_i u_{i\bar i},\qquad
 J\partial_{x_j}=\partial_{y_j},\quad
 J\partial_{y_j}=-\partial_{x_j}.
\]
Let $B_r(a)$ be the Euclidean ball of radius $r$ centered at $a$,
and write $B_r=B_r(0)$ and $d(z)=\operatorname{dist}(z,\partial\Omega)$.
The outward unit normal is denoted by $\nu$. Integrals over $\Omega$
use real $2n$-dimensional Lebesgue measure.
All H\"older and Sobolev norms are Euclidean unless another metric is
specified. When the domain is omitted, $L^r=L^r(\Omega)$ for
$1\le r\le\infty$. On approximating domains, $L^r=L^r(\Omega_j)$.
For a continuous tensor field $T$ on a fixed compact set $K$, write
\[
 \omega_T(t)=\sup_{\substack{x,y\in K\\|x-y|\le t}}|T(x)-T(y)|.
\] A constant $C$ may change from line to line. Unless stated otherwise,
it depends only on $n$, $\alpha$, $\lambda$, $\Lambda$, $\|f\|_{C^\alpha}$,
a fixed finite $C^3$ boundary atlas, a positive lower bound for the
Levi form on unit complex tangential vectors, and $\|\Phi\|_{C^3}$
for the boundary extension $\Phi$.
The constants may also depend on fixed continuity moduli for the third
derivatives of the boundary charts and boundary extension. Dependence on these continuity
moduli does not require a Dini condition or a positive H\"older exponent
for the third derivatives. 

The a priori estimates are first proved for smooth data on strongly
pseudoconvex domains, where smooth solutions exist by \cite{CKNS}.
The constants are uniform under smooth approximation of the fixed $C^3$ data.

\subsection{Defining functions and smooth approximation}

We will construct barriers, prove a weighted integral estimate, and
approximate the domain and data by smooth ones. The barrier estimates follow from
\cite[Section~3.2, equation~(3.8)]{Guan2009} and
\cite[Lemma~2.1]{Guan}, together with the arithmetic--geometric mean
inequality.

\begin{lemma}\label{lem:geometric-barriers}

Let $\Omega\Subset\mathbb C^n$ have strongly pseudoconvex $C^3$
boundary, and let $\varphi\in C^3(\partial\Omega)$ be real-valued.
There are real-valued $C^3$ functions $\Phi$ and $\rho$ on a fixed open
neighborhood of $\bar\Omega$ such that
\[
 \begin{gathered}
  \Phi|_{\partial\Omega}=\varphi,\qquad
  \rho|_{\partial\Omega}=0,\qquad \rho<0\quad\text{in }\Omega,\\
  (\rho_{i\bar j})\ge c_\rho I,\qquad
  |D\rho|\ge c_\rho>0\quad\text{on }\partial\Omega.
 \end{gathered}
\]
There are fixed constants $A,B,c,C>0$ such that every smooth strictly
plurisubharmonic solution $u$ of \eqref{intro:dirichlet} with
$\lambda\le f\le\Lambda$ satisfies
\begin{equation}
 \Phi+A\rho\le u\le\Phi-B\rho.                              \label{eq:CG4}
\end{equation}
Moreover,
\begin{equation}
 \|u\|_{C^0(\overline\Omega)}+
 \|Du\|_{C^0(\partial\Omega)}\le C.                         \label{eq:CG5}
\end{equation}
The function $\underline u=\Phi+A\rho$ is strictly plurisubharmonic, and
\begin{equation}
 L_u(u-\underline u)=n-\sum_{i,j}u^{i\bar j}\underline u_{i\bar j}
       \le-c(1+\sum_i u^{i\bar i}).                                       \label{eq:CG6}
\end{equation}
The constants depend only on the fixed geometric and boundary data
and the bounds $\lambda$ and $\Lambda$ for the right-hand side.
\end{lemma}

\begin{proof}
By the boundary extension construction in
\cite[Lemma~6.38]{GT}, we obtain a $C^3$ function $\Phi$ on a fixed
neighborhood of $\bar\Omega$ such that
\begin{equation}
 \Phi|_{\partial\Omega}=\varphi,
 \qquad \|\Phi\|_{C^3}\le C.                         \label{eq:CG1}
\end{equation}
Let $\rho_0$ be a real-valued $C^3$ defining function on a fixed
neighborhood of $\partial\Omega$. In this neighborhood, $\rho_0$ is
negative in $\Omega$ and positive outside $\bar\Omega$, with
$\rho_0=0$ and $|D\rho_0|>0$ on $\partial\Omega$.
By strong pseudoconvexity, $\rho_0+K\rho_0^2$ is uniformly strictly
plurisubharmonic on a fixed boundary collar for sufficiently large $K$. Let $R\ge1$ be fixed so that
$B_R$ contains a neighborhood of $\bar\Omega$. For a fixed sufficiently
small $a>0$ with $Ka<1/4$, set
$q(z)=\frac{a}{16R^2}|z|^2-\frac{3a}{8}$.
Then
\[
 \begin{aligned}
 -\frac{3a}{8}\le q&\le-\frac{5a}{16}<0
                    &&\text{in }B_R,\\
 q&<\rho_0+K\rho_0^2=0
                    &&\text{on }\{\rho_0=0\},\\
 \rho_0+K\rho_0^2&\le-\frac{3a}{4}<-\frac{3a}{8}\le q
                    &&\text{on }\{\rho_0=-a\}.
 \end{aligned}
\]
Let
$\operatorname{Max}_{a/32}(s,t)$ be a fixed smooth convex function,
nondecreasing in $s$ and $t$, such that
\[
 \begin{gathered}
 \partial_s\operatorname{Max}_{a/32}
       +\partial_t\operatorname{Max}_{a/32}=1,\\
 \operatorname{Max}_{a/32}(s,t)=\max\{s,t\}
                      \quad\text{if }|s-t|\ge a/32.
 \end{gathered}
\]
Set
\[
 \rho=
 \begin{cases}
  \operatorname{Max}_{a/32}(\rho_0+K\rho_0^2,q),
                                      &-a<\rho_0\le0,\\
  q,                                  &\text{farther inside }\Omega.
 \end{cases}
\]
Near $\{\rho_0=-a\}$, we have $\rho=q$. The formula
$\rho=\rho_0+K\rho_0^2$ near $\partial\Omega$ also defines $\rho$
across the boundary. By the convexity and monotonicity of
$\operatorname{Max}_{a/32}$, we obtain a $C^3$ defining function
with bounded $C^3$ norm satisfying
\begin{equation}
 (\rho_{i\bar j})\ge c_\rho I\quad\hbox{near }\bar\Omega,
 \qquad |D\rho|\ge c_\rho\quad\hbox{on }\partial\Omega.
                                                        \label{eq:CG2}
\end{equation}
We fix $A,B>0$ sufficiently large so that
\eqref{eq:CG1}--\eqref{eq:CG2} give
\begin{equation}
  (\underline u_{i\bar j})\ge c_0I,\qquad
 \det(\underline u_{i\bar j})\ge3^n\Lambda,           \label{eq:CG3}
\end{equation}
and $\Delta(\Phi-B\rho)<0$. Using the comparison principle
\cite[Theorem~A]{BT} and the maximum principle
\cite[Theorem~3.1]{GT}, we obtain \eqref{eq:CG4} and the $C^0$ bound.
Differentiating \eqref{eq:CG4} at the boundary gives a bound for $u_\nu$.
Since the tangential derivatives are determined by $\varphi$,
we obtain \eqref{eq:CG5}.

By \eqref{eq:CG3}, $f\le\Lambda$ and the arithmetic--geometric mean inequality, we have
\[
 L_u\underline u\ge3n;\qquad
 L_u\underline u\ge c_0\sum_i u^{i\bar i};\qquad L_uu=n.
\]
It follows that
$L_u(u-\underline u)=n-L_u\underline u
\le-\frac n2-\frac{c_0}{2}\sum_i u^{i\bar i}$,
which proves \eqref{eq:CG6}.
\end{proof}

\Needspace{5\baselineskip}
For $\psi\in C^1(\bar\Omega)$ with $\psi=\Phi$ on $\partial\Omega$, write
$\gamma=(\psi-\Phi)/\rho$, extended continuously to the boundary by
$\gamma=(\psi_\nu-\Phi_\nu)/\rho_\nu$.

\begin{lemma}
\label{lem:boundary-quotient-base}
We use $d^c=\frac{i}{2}(\bar\partial-\partial)$, so that
$dd^c=i\partial\bar\partial$, and set $\omega_0=dd^c|z|^2$.
Let $\psi\in C^\infty(\bar\Omega)$ be plurisubharmonic with
$\psi=\Phi$ on $\partial\Omega$, and suppose
$|\gamma|\le B$.
Put $\omega_\psi=dd^c\psi$.
Then, for every integer $0\le k\le n-1$,
\begin{equation}
 \int_\Omega(-\rho)\,i\partial\gamma\wedge\bar\partial\gamma
       \wedge\omega_\psi^k\wedge\omega_0^{n-k-1}\le C.
 \label{eq:boundary-quotient-base}
\end{equation}
The constant depends only on $B$ and the fixed geometric and boundary
bounds. In particular, the cases $k=0$ and $k=n-1$ give
\[
 \int_\Omega(-\rho)|D\gamma|^2\le C,
 \qquad \int_\Omega(-\rho)\,Q|\partial\gamma|_g^2\le C,
 \quad g=(\psi_{i\bar j}),\quad Q=\det g.
\]%
The second estimate applies when $g>0$, with
$|\partial\gamma|_g^2=g^{i\bar j}\gamma_i\gamma_{\bar j}$.
\end{lemma}
\begin{proof}
Let $A>B$ be fixed and sufficiently large that
$v=\Phi+A\rho$ satisfies $dd^cv\ge c\omega_0$.
Then $v\le\psi$ and $v=\psi$ on $\partial\Omega$, so
$\partial_\nu(\psi-v)\le0$. By Stokes' theorem,
\[
 \begin{aligned}
 \frac{d}{ds}\int_\Omega[dd^c(2v+s(\psi-v))]^n
 &=n\int_{\partial\Omega}d^c(\psi-v)\wedge
       [dd^c(2v+s(\psi-v))]^{n-1}\\
 &\le0,\qquad0\le s\le1.
 \end{aligned}
\]
Integrating in $s$, we obtain
\[
 \int_\Omega(\omega_\psi+dd^cv)^n
       \le\int_\Omega(2dd^cv)^n\le C.
\]
Positivity of the mixed terms and $dd^cv\ge c\omega_0$ give
\[
 \int_\Omega\omega_\psi^j\wedge\omega_0^{n-j}\le C,
 \qquad 0\le j\le n.
\]
Set $T_k=\omega_\psi^k\wedge\omega_0^{n-k-1}$.
This form is closed and positive. By the product rule, we have
\[
 2(-\rho)\,i\partial\gamma\wedge\bar\partial\gamma
   =-dd^c(\rho\gamma^2)+2\gamma\,dd^c(\psi-\Phi)-\gamma^2dd^c\rho.
\]%
Stokes' theorem \cite[Theorem~10.50]{Rudin} then gives
\begin{align*}
 2\int_\Omega(-\rho)\,i\partial\gamma\wedge\bar\partial\gamma\wedge T_k
 ={}&-\int_{\partial\Omega}\gamma^2d^c\rho\wedge T_k
  +2\int_\Omega\gamma(\omega_\psi-dd^c\Phi)\wedge T_k\\
 &-\int_\Omega\gamma^2dd^c\rho\wedge T_k.
\end{align*}
The boundary term is nonpositive because
$d^c\rho\wedge T_k$ is nonnegative with respect to the induced
boundary orientation.
Using the mixed Monge--Amp\`ere mass bounds, $|\gamma|\le B$ and the fixed
$C^2$ bounds for $\Phi,\rho$, we bound both interior terms by $C$.
Thus \eqref{eq:boundary-quotient-base} follows.
\end{proof}

Using the extension results in \cite[Lemmas~6.37 and~6.38]{GT}
and mollification \cite[Section~7.2]{GT}, we approximate the $C^3$
data while preserving the geometric bounds and the H\"older bound
of the right-hand side.
For related results on variable domains, we refer to
Gong--Kim \cite[Section~2]{GK}.

\begin{lemma}\label{lem:smooth-data}
Let $\rho$ and $\Phi$ be as in Lemma~\ref{lem:geometric-barriers},
and suppose that $f\in C^\alpha(\bar\Omega)$, $0<\alpha<1$, with
$0<\lambda\le f\le\Lambda$. There are a fixed bounded open
neighborhood $U$ of $\bar\Omega$ and smooth functions
$\rho_j,\Phi_j,f_j$ on a neighborhood of $\bar U$ with the following properties.
\begin{enumerate}
\item We have
$\|\rho_j-\rho\|_{C^3(\bar U)}
   +\|\Phi_j-\Phi\|_{C^3(\bar U)}\longrightarrow0$.
\item\label{item:smooth-domains} For a fixed $z_0\in\Omega$, let $\Omega_j$
be the connected component of $\{\rho_j<0\}\cap U$ containing $z_0$,
and set $\varphi_j=\Phi_j|_{\partial\Omega_j}$.
Each $\Omega_j$ has smooth strongly pseudoconvex boundary.
\item The boundary values $\varphi_j$ converge to $\varphi$ in $C^3$
when expressed in the same local tangential coordinates. The functions $\rho_j$ and
$\Phi_j$ have uniformly bounded $C^3$ norms and common
continuity moduli for their third derivatives.
The lower bounds in \eqref{eq:CG2} are independent of $j$.
\item The functions $f_j$ satisfy
\[
 \lambda\le f_j\le\Lambda\quad\text{on }\bar U,\qquad
 [f_j]_{C^\alpha(\bar U)}\le[f]_{C^\alpha(\bar\Omega)},\qquad
 \|f_j-f\|_{C^0(\bar\Omega)}\longrightarrow0.
\]
\end{enumerate}
\end{lemma}

\begin{proof}
Let $\eta_\varepsilon$ be a standard nonnegative mollifier supported
in $B_\varepsilon$. On a fixed bounded neighborhood $U$ of
$\bar\Omega$, with $\bar U$ contained in the common extension
neighborhood and $\rho>0$ on $\partial U$, set
$\rho_j=\rho*\eta_{\varepsilon_j}$ and
$\Phi_j=\Phi*\eta_{\varepsilon_j}$, where $\varepsilon_j\downarrow0$.
By mollification \cite[Section~7.2]{GT}, we have $C^3$ convergence.
For sufficiently small $\varepsilon_j$, the implicit function theorem
\cite[Theorem~9.28]{Rudin} gives smooth $\partial\Omega_j$ whose
graph functions converge in $C^3$ on fixed smaller charts.
Composing with $\Phi_j$, we obtain $\varphi_j\to\varphi$ in $C^3$
on these charts. The common $C^3$ bounds
and third-derivative moduli follow from the convolution formulas and
the chain rule.
The positive lower bounds in \eqref{eq:CG2} also persist.

Let
\[
 \widehat f(z)=\min\left\{\Lambda,
   \inf_{x\in\bar\Omega}
   \{f(x)+[f]_{C^\alpha(\bar\Omega)}|z-x|^\alpha\}\right\},
 \qquad f_j=\widehat f*\eta_{\varepsilon_j}.
\]
Since $|x-y|^\alpha$ is a metric, we have $\widehat f=f$ on
$\bar\Omega$, $\lambda\le\widehat f\le\Lambda$, and
$[\widehat f]_{C^\alpha(\mathbb C^n)}\le[f]_{C^\alpha(\bar\Omega)}$.
These bounds are preserved by nonnegative convolution. Moreover,
$\|f_j-f\|_{C^0(\bar\Omega)}
 \le[f]_{C^\alpha(\bar\Omega)}\varepsilon_j^\alpha\to0$.
This proves the lemma.
\end{proof}

\subsection{Elliptic estimates and comparison}\label{sec:uniform-elliptic}

We will derive the global $C^{2,\alpha}$ estimate from a Hessian
bound and establish a stability estimate with respect to the data.
The regularity estimate follows from \cite[Theorem~1.4(b)]{SS}
and the linear Schauder estimates \cite[Sections~6.1--6.2]{GT}.

\begin{proposition}\label{prop:elliptic-finish}
Let $0<\alpha<1$ and $0<\lambda\le\Lambda$.
On a fixed smooth bounded domain $\Omega$,
let $u\in C^\infty(\bar\Omega)$ be a strictly plurisubharmonic
solution of \eqref{intro:dirichlet} with fixed smooth real-valued
boundary data $\varphi$. Suppose
\begin{equation}
 \|u\|_{C^1(\bar\Omega)}+\|D^2u\|_{L^\infty}\le M,
 \qquad \lambda\le\det {\bigl(u_{i\bar j}\bigr)}=f\le\Lambda.
 \label{eq:FEU1}
\end{equation}
Then
\begin{equation}
                    \|u\|_{C^{2,\alpha}(\bar\Omega)}\le C,
 \label{eq:FEU2}
\end{equation}
where $C$ depends on $M,\lambda,\Lambda,\|f\|_{C^\alpha(\bar\Omega)},\alpha$,
the fixed domain and the fixed boundary data, but not on higher derivatives of $f$.

\end{proposition}

\begin{proof}
For a real symmetric matrix $A$ in coordinates $(x,y)$, set
\[
 \mathcal P(A)_{i\bar j}
 =\tfrac14(A_{x_ix_j}+A_{y_iy_j}
                +iA_{x_iy_j}-iA_{y_ix_j}).
\]
We have $\mathcal P(D^2u)=(u_{i\bar j})$. By \eqref{eq:FEU1},
$aI\le\mathcal P(D^2u)\le bI$ for fixed $0<a<b$.
For $F(A)=\log\det\mathcal P(A)$, we use the extension
construction in \cite[Section~2]{YW} and define
\begin{equation}
 \widehat F(A)=\inf_{aI\le B\le bI}
       \{\log\det B-n+\operatorname{tr}(B^{-1}\mathcal P(A))\},
 \label{eq:FEU3a}
\end{equation}
where $B$ ranges over Hermitian matrices. The operator $\widehat F$ is
concave, with ellipticity constants
$1/(4b)$ and $1/(4a)$. Thus $\widehat F(D^2u)=F(D^2u)=\log f$.
We apply \cite[Theorem~1.4(b)]{SS} to $-u$ with the convex operator
$A\mapsto-\widehat F(-A)+\widehat F(0)$, which vanishes at zero,
and obtain
\begin{equation}
 \|u\|_{C^{2,\gamma}(\bar\Omega)}\le C_\gamma       \label{eq:FEU4}
\end{equation}
for some fixed $0<\gamma<\alpha$.

For sufficiently small $s>0$, we use boundary neighborhoods
$\Omega\cap B_{2s}(p)$ with $p\in\partial\Omega$ and interior balls
$B_{2s}(p)\Subset\Omega$. For $x_0\in\Omega\cap B_{2s}(p)$,
set $A_0=D^2u(x_0)$ and define
\[
 R_{A_0}(A)=F(A)-F(A_0)
       -\operatorname{tr}\bigl(\mathcal P(A_0)^{-1}\mathcal P(A-A_0)\bigr).
\]
By Taylor's formula \cite[Theorem~5.15]{Rudin} and
\eqref{eq:FEU4}, we have
\[
 \begin{gathered}
  [R_{A_0}(D^2u)]_{C^\alpha(\Omega\cap B_{2s}(p))}
       \le Cs^\gamma[D^2u]_{C^\alpha(\Omega\cap B_{2s}(p))},\\
  u^{i\bar j}(x_0)u_{i\bar j}
       =\log f-F(A_0)+n-R_{A_0}(D^2u).
 \end{gathered}
\]
Using the local Schauder estimates \cite[Sections~6.1--6.2]{GT},
we get
\begin{equation}
 [D^2u]_{C^\alpha(\Omega\cap B_s(p))}
       \le Cs^\gamma[D^2u]_{C^\alpha(\Omega\cap B_{2s}(p))}+C_s.
 \label{eq:FEU-local-schauder}
\end{equation}
Here $C$ depends only on the stated data, and $C_s$ may also depend
on $s$.
For $|x-y|<s/8$, we use \eqref{eq:FEU-local-schauder}
on $\Omega\cap B_s(p)$ if $d(x)<s/2$, where $p$ is a nearest boundary
point to $x$. If $d(x)\ge s/2$, we use the same estimate on
$B_{s/8}(x)$, since $B_{s/4}(x)\Subset\Omega$. For $|x-y|\ge s/8$,
the H\"older quotient is at most $2M(8/s)^\alpha$. Taking the
supremum gives
\begin{equation}
 [D^2u]_{C^\alpha(\bar\Omega)}
       \le Cs^\gamma[D^2u]_{C^\alpha(\bar\Omega)}+C_s. \label{eq:FEU9}
\end{equation}
Taking $Cs^\gamma\le1/2$, we absorb the first term on the right
of \eqref{eq:FEU9} and obtain \eqref{eq:FEU2}.
\end{proof}

Linearizing between two solutions and applying the Schauder
estimate gives the following stability bound.

\begin{proposition}\label{prop:holder-stability}
Let $D\Subset\mathbb C^n$ have $C^{2,\alpha}$ boundary,
$0<\alpha<1$, and fix $\kappa,M>0$.
Suppose that real-valued functions $u,v\in C^{2,\alpha}(\bar D)$ satisfy
\begin{equation}
 \bigl(u_{i\bar j}\bigr),\bigl(v_{i\bar j}\bigr)\ge\kappa I,\qquad
 \|u\|_{C^{2,\alpha}(\bar D)}+
                 \|v\|_{C^{2,\alpha}(\bar D)}\le M.
 \label{eq:stability-hessian-data}
\end{equation}
Set $f=\det {\bigl(u_{i\bar j}\bigr)}$ and $g=\det {\bigl(v_{i\bar j}\bigr)}$. Then
\begin{equation}
 \|u-v\|_{C^{2,\alpha}(\bar D)}
 \le C\left(\|f-g\|_{C^\alpha(\bar D)}
          +\|u-v\|_{C^{2,\alpha}(\partial D)}\right),
                                                               \label{eq:CO1}
\end{equation}
where $C$ depends only on $D,n,\alpha,\kappa,M$.
\end{proposition}

\begin{proof}
Let $w=u-v$ and $H_t=((1-t)v_{i\bar j}+tu_{i\bar j})$.
Linearizing the determinant, we obtain
\begin{equation}
 B^{i\bar j}w_{i\bar j}=f-g,\qquad
 B^{i\bar j}=\int_0^1\det H_t\,(H_t^{-1})_{ji}\,dt.
 \label{eq:CO2}
\end{equation}
By \eqref{eq:stability-hessian-data}, the associated real operator
is uniformly elliptic and its coefficients have uniformly bounded
$C^\alpha$ norms. The maximum principle and the global Schauder
estimate \cite[Theorem~6.6]{GT} give \eqref{eq:CO1}.
\end{proof}

\subsection{Complex tangential Hessian estimates}\label{sec:tangential-block}

We will prove the complex tangential Hessian lower bound in
\cite[Proposition~2.2]{Guan} and \cite[Proposition~3.9]{Guan2009}
under $C^3$ assumptions on the boundary and boundary values.
The proof in \cite{Guan2009} uses a $C^4$ defining function.
Here a $C^1$ coefficient in the comparison function is regularized,
and the errors are controlled by the continuity moduli of third
derivatives.

\begin{proposition}\label{prop:tangential-block}
Let $\Omega\Subset\mathbb C^n$, $n\ge2$, have strongly pseudoconvex
$C^3$ boundary. Choose a fixed neighborhood $N$ of $\partial\Omega$,
and let $\rho,\Phi\in C^3(N)$ be real-valued. Assume that $\rho$ is a
defining function with $\rho<0$ on $N\cap\Omega$,
$\bigl(\rho_{i\bar j}\bigr)\ge c_\rho I$, and $|D\rho|\ge c_\rho>0$
on $N$. For a fixed sufficiently small $\delta>0$, let
$\mathcal C_\delta=\{z\in N\cap\Omega:0<-\rho(z)<\delta\}$, with
$\overline{\mathcal C_\delta}\subset N$. Suppose that a real-valued
$u\in C^2(\overline{\mathcal C_\delta})$ satisfies
\begin{equation}
 \begin{gathered}
  {\bigl(u_{i\bar j}\bigr)}>0,\qquad \det {\bigl(u_{i\bar j}\bigr)}\ge\lambda>0\quad\hbox{in }\mathcal C_\delta,\\
  u=\Phi\quad\hbox{on }\partial\Omega,\qquad
  |u-\Phi|\le B(-\rho)\quad\hbox{in }\mathcal C_\delta.
 \end{gathered}                                                \label{eq:GC1}
\end{equation}
Write
\[
 T_p^{1,0}\partial\Omega
   =\{\xi\in\mathbb C^n:\textstyle\sum_i\rho_i(p)\xi_i=0\},\qquad
 (v_{i\bar j})(\xi,\bar\eta)
   =\sum_{i,j}v_{i\bar j}\xi_i\overline{\eta_j}
\]
for real-valued $v\in C^2$.
Then there is a constant $c_T>0$ such that
\[
 (u_{i\bar j}(p))(\xi,\bar\xi)\ge c_T|\xi|^2
 \qquad(p\in\partial\Omega,\ \xi\in T_p^{1,0}\partial\Omega).
\]
The constant $c_T$ depends only on $n$, $\delta$, $c_\rho$, the fixed
$C^3$ geometric and boundary data, their third-derivative continuity
moduli, $B$, and $\lambda$. The estimate is uniform for smooth
approximations with these bounds and moduli.
\end{proposition}

\begin{proof}
Fix $p\in\partial\Omega$ and a unit vector
$\xi_p\in T_p^{1,0}\partial\Omega$. On a ball of fixed radius about
$p$, project $\xi_p$ onto $\ker\partial\rho$ and normalize:
\[
 \widetilde\xi_i(z)=\xi_{p,i}
 -\frac{\sum_k\rho_k(z)\xi_{p,k}}{|\partial\rho(z)|^2}
       \overline{\rho_i(z)},\qquad
 \xi=\widetilde\xi/|\widetilde\xi|.
\]
The radius can be chosen uniformly so that $|\widetilde\xi|\ge1/2$.
Thus $\xi$ is a $C^2$ unit complex tangential vector field with
uniform bounds and $\xi(p)=\xi_p$. Set
\begin{equation}
 a_*=-\frac{(\Phi_{i\bar j})(\xi,\bar\xi)}
                   {(\rho_{i\bar j})(\xi,\bar\xi)}.
 \label{eq:GC2}
\end{equation}
We have $\|a_*\|_{C^1}\le C$ and
\[
 \omega_{Da_*}(t)\le
 C\{\omega_{D^3\Phi}(Ct)+\omega_{D^3\rho}(Ct)+t\}.
\]
Choose a fixed $L<0$ with $\bigl(\Phi_{i\bar j}+L\rho_{i\bar j}\bigr)\le-I$, and let
$0\le\zeta\le1$ be supported in this ball and equal to one near $p$.
We extend $a=\zeta a_*+(1-\zeta)L$ by $L$ outside the ball.
Since $L\le a_*$, we obtain $a(p)=a_*(p)$ and
\begin{equation}
 (\Phi_{i\bar j}+a\rho_{i\bar j})(\xi_z,\bar\xi_z)\le0
 \quad\text{for some unit }\xi_z\in\ker\partial\rho(z)
 \text{ at each }z\in\mathcal C_\delta.
 \label{eq:GC3}
\end{equation}
The $C^1$ bound for $a$ and a continuity modulus $\omega$ for $Da$
are uniform in $(p,\xi_p)$, and $\omega(t)\to0$ as $t\downarrow0$.

For a fixed smooth nonnegative even mollifier of integral one, put
$a_\ell=a*\eta_\ell$. Taylor's formula and the cancellation of
$D\eta_\ell$ give
\begin{equation}
 \|a_\ell-a\|_{L^\infty(\mathbb C^n)}\le C\ell\omega(\ell),\qquad
 \|Da_\ell\|_{L^\infty(\mathbb C^n)}\le C,\qquad
 \|D^2a_\ell\|_{L^\infty(\mathbb C^n)}\le C\omega(\ell)/\ell.
 \label{eq:GC4}
\end{equation}
Choose fixed $M_0,K_0$ sufficiently large, with
$M_0\delta\ge B+\|a\|_{L^\infty(\mathbb C^n)}+1$, and set
\begin{equation}
 M_\ell=M_0+K_0\omega(\ell)/\ell,\qquad
 \widehat{\Phi}=\Phi+(a_\ell+\ell)\rho+M_\ell\rho^2,
 \qquad0<\ell\le1.
 \label{eq:GC9}
\end{equation}
Differentiation yields
\begin{equation}
 \begin{aligned}
 \widehat{\Phi}_{i\bar j}={}&\Phi_{i\bar j}+(a_\ell+\ell)\rho_{i\bar j}
       +2M_\ell\rho\,\rho_{i\bar j}+\rho(a_\ell)_{i\bar j}\\
 &+(a_\ell)_i\rho_{\bar j}+\rho_i(a_\ell)_{\bar j}
       +2M_\ell\rho_i\rho_{\bar j}.
 \end{aligned}
 \label{eq:GC11}
\end{equation}
By \eqref{eq:GC4}, a fixed sufficiently large $K_0$ ensures
\begin{equation}
 2M_\ell\rho(\rho_{i\bar j})+\rho((a_\ell)_{i\bar j})
 \le (-\rho)\{C\omega(\ell)/\ell-2c_\rho M_\ell\}I\le0.
 \label{eq:GC12}
\end{equation}
The last three terms in \eqref{eq:GC11} vanish on complex tangent
vectors. Hence \eqref{eq:GC3}--\eqref{eq:GC12} give
\[
 \begin{gathered}
 (\widehat{\Phi}_{i\bar j})(\xi_z,\bar\xi_z)\le C\ell(1+\omega(\ell)),\qquad
 (\widehat{\Phi}_{i\bar j})(\tau,\bar\tau)\le C
 \quad(\tau\in\ker\partial\rho,\ |\tau|=1),\\
 (\widehat{\Phi}_{i\bar j})(\nu_z,\bar\nu_z)\le C(1+M_\ell),
 \end{gathered}
\]
where $\nu_z$ is a unit complex normal vector. At points where
$(\widehat{\Phi}_{i\bar j})>0$, Hadamard's inequality \cite[Section~7.8]{HJ} therefore gives
\begin{equation}
 \begin{aligned}
 \det(\widehat{\Phi}_{i\bar j})
 &\le C\ell(1+\omega(\ell))(1+M_\ell)\\
 &=C(1+\omega(\ell))\{(1+M_0)\ell+K_0\omega(\ell)\}
 \xrightarrow{\ell\downarrow0}0.
 \end{aligned}
 \label{eq:GC16}
\end{equation}
Fix $\ell>0$ sufficiently small that
\[
 \|a_\ell-a\|_{L^\infty(\mathbb C^n)}\le\ell/2,\qquad
 C(1+\omega(\ell))\{(1+M_0)\ell+K_0\omega(\ell)\}
 \le\lambda/2.
\]%
On the inner boundary, \eqref{eq:GC1} and the choice of $M_0$ give
\[
 \widehat{\Phi}-u\ge M_0\delta^2-(B+\|a\|_{L^\infty(\mathbb C^n)}+1)\delta\ge0
 \quad\text{on }\{\rho=-\delta\}.
\]
Also $\widehat{\Phi}=u$ on $\partial\Omega$. A positive interior maximum of
$u-\widehat{\Phi}$ would imply $(\widehat{\Phi}_{i\bar j})\ge(u_{i\bar j})>0$ and
$\det(\widehat{\Phi}_{i\bar j})\ge\lambda$, contradicting \eqref{eq:GC16}.
Thus $u\le \widehat{\Phi}$ throughout the collar.

Normal differentiation of $u\le \widehat{\Phi}$ at $p$ gives
$\gamma(p)\ge a_\ell(p)+\ell$. Differentiating $u-\Phi=0$ twice
in complex tangential directions, we obtain
\begin{equation}
 \begin{aligned}
 (u_{i\bar j}(p))(\xi_p,\bar\xi_p)
 &=(\gamma(p)-a_*(p))(\rho_{i\bar j}(p))(\xi_p,\bar\xi_p)\\
 &\ge\frac\ell2(\rho_{i\bar j}(p))(\xi_p,\bar\xi_p).
 \end{aligned}
 \label{eq:GC19}
\end{equation}
Consequently
\begin{equation}
 (u_{i\bar j})|_{T^{1,0}\partial\Omega}\ge c_T I,
 \qquad c_T=c_\rho\ell/2>0.
 \label{eq:GC20}
\end{equation}
All choices are uniform in $(p,\xi_p)$ and for smooth approximations
with the stated common bounds and continuity moduli.
\end{proof}

For zero boundary values, strong pseudoconvexity is necessary for a
classical PSH solution with positive determinant on $\bar\Omega$,
as observed in \cite[p.~689]{Guan}.
Indeed, if $u\in\PSH(\Omega)\cap C^2(\bar\Omega)$, $u=0$ on
$\partial\Omega$, and $\det(u_{i\bar j})\ge\lambda>0$, then
$(u_{i\bar j})>0$ on $\bar\Omega$. The maximum principle and Hopf
lemma \cite[Lemma~3.4]{GT} give $u<0$ in $\Omega$ and $u_\nu>0$
on $\partial\Omega$. For a defining function $\rho<0$, the boundary
identity is
\[
 (u_{i\bar j})(\xi,\bar\xi)
 =\frac{u_\nu}{\rho_\nu}(\rho_{i\bar j})(\xi,\bar\xi),
 \qquad \xi\in T^{1,0}\partial\Omega.
\]
Hence the Levi form is positive definite.

\subsection{Uniform estimates and compactness}\label{sec:analytic-uniformity}

Using the approximations of Lemma~\ref{lem:smooth-data}, we will
establish uniform analytic estimates and show that a common
$C^{2,\alpha}$ bound passes to the limit.

Let $P_t$ be the Neumann heat semigroup on $\Omega_j$.
For $h\in C(\bar\Omega_j)$, the function $v(x,t)=(P_th)(x)$ solves
\begin{equation}
 \begin{cases}
  \partial_t v=\Delta v&\text{in }\Omega_j\times(0,\infty),\\
  \partial_\nu v=0&\text{on }\partial\Omega_j\times(0,\infty),\\
  v(\cdot,0)=h&\text{on }\bar\Omega_j,
 \end{cases}
 \label{eq:neumann-heat-definition}
\end{equation}
where $\Delta$ is the real Euclidean Laplacian.
The corresponding Neumann heat kernel $p_t(x,y)$ satisfies
\[
 (P_th)(x)=\int_{\Omega_j}p_t(x,y)h(y)\,dy,
       \qquad x\in\Omega_j,\quad t>0.
\]
The same formula defines $P_t$ on $L^r(\Omega_j)$ for
$1\le r\le\infty$.
The common $C^2$ boundary charts give uniform bounds for the
second fundamental form and uniform lower bounds for the volumes
of small intrinsic balls, as required in \cite[Section~2]{WY}. On the fixed boundary neighborhoods, we choose
functions with inward normal derivative one on $\partial\Omega_j$
and uniformly bounded $C^2$ norms. By It\^o's formula for reflecting Brownian motion, we obtain a uniform
bound for the exponential moment of boundary local time in
\cite[Eq.~(2.3)]{WY}.
Using \cite[equation~(2.2)]{WY} and the $L^1$ contraction property
of $P_t$, we obtain
$\|e^{-t}P_t\|_{L^1\to L^\infty}\le Ct^{-n}$ for $t>0$.
For $\kappa\in\mathbb R$ and a real-valued Lipschitz function $\psi$
with Lipschitz constant at most one,
the energy identity for the Neumann heat equation gives
\(\|e^{\kappa\psi}P_te^{-\kappa\psi}\|_{L^2\to L^2}\le e^{\kappa^2t}\).
Optimizing in $\kappa$, we obtain the Davies--Gaffney estimates
\cite[Section~3, Eq.~(3.2)]{CS} for $P_t$ and $e^{-t}P_t$,
uniformly in $j$.
By \cite[Theorem~4.1]{CS}, applied to $e^{-t}P_t$, and
\cite[Theorem~1.1 and Eq.~(1.2)]{WY}, we have, for $0<t\le1$,
\begin{equation}
 p_t(x,y)\le Ct^{-n}e^{-|x-y|^2/(Ct)},\qquad
 |\nabla_xp_t(x,y)|+|\nabla_yp_t(x,y)|
     \le Ct^{-(2n+1)/2}e^{-|x-y|^2/(Ct)}.
 \label{eq:uniform-neumann-kernel}
\end{equation}
The intrinsic distance is no smaller than the Euclidean distance.
The derivative bound in the second variable follows by symmetry.

In a boundary chart $x_{2n}=\psi_j(y)$, the map
$(y,r)\mapsto(y,r+\psi_j(y))$ has Jacobian determinant one and uniform
bi-Lipschitz bounds. By the chain rule, the first-derivative
bounds in \eqref{eq:uniform-neumann-kernel} remain uniform in the graph
coordinates, since $D\psi_j$ and the cutoff derivatives are uniformly
bounded. For $f\in C^\alpha(\bar\Omega_j)$ and $0<\alpha<1$,
using \eqref{eq:uniform-neumann-kernel} and the identities
$\int_{\Omega_j} p_t(x,y)\,dy=1$ and
$\int_{\Omega_j}\nabla_xp_t(x,y)\,dy=0$, we get
\[
 \|P_tf-f\|_{L^\infty}\le Ct^{\alpha/2}\|f\|_{C^\alpha},\qquad
 \|\nabla P_tf\|_{L^\infty}\le
                 Ct^{(\alpha-1)/2}\|f\|_{C^\alpha}.
\]
Nearby points in each boundary chart can be joined inside the domain
by paths of length at most $C|x-y|$. Integrating $\nabla P_tf$
along these paths gives
\[
 |P_tf(x)-P_tf(y)|
 \le Ct^{(\alpha-1)/2}\|f\|_{C^\alpha}|x-y|
 \quad\text{if }|x-y|\le\sqrt t.
\]
For $|x-y|\ge\sqrt t$, the bound for $\|P_tf-f\|_{L^\infty(\Omega_j)}$ gives
the same H\"older estimate. Hence
$\|P_tf\|_{C^\alpha}\le C\|f\|_{C^\alpha}$.

Using the semigroup property and integration by parts with the
Neumann boundary condition, we have
\[
 \partial_tP_tf(x)
 =-\int_{\Omega_j}\nabla_y p_{t/2}(x,y)\cdot
                         \nabla P_{t/2}f(y)\,dy.
\]
Combining the gradient estimate for the Neumann heat kernel with the
bound for $\nabla P_tf$ above and using the semigroup property once more,
we obtain
\[
 \|\partial_tP_tf\|_{L^\infty}
 +t^{1/2}\|\nabla\partial_tP_tf\|_{L^\infty}
 \le Ct^{-1+\alpha/2}[f]_{C^\alpha}.
\]
Integrating the Gaussian bound also gives
$\|P_tf-f(x_0)\|_{L^\infty(\Omega_j\cap B_{C\sqrt t}(x_0))}
 \le Ct^{\alpha/2}[f]_{C^\alpha}$.
By the scaled local Neumann Schauder estimate
\cite[Lemma~6.29 and Theorem~6.30]{GT} with exponent $1/2$, we get
\begin{equation}
 \|D^2P_tf\|_{L^\infty}+\|\partial_tP_tf\|_{L^\infty}
 \le Ct^{-1+\alpha/2}[f]_{C^\alpha},\qquad 0<t\le1.
 \label{eq:neumann-second-derivative}
\end{equation}
Indeed, $\partial_tP_tf$ has scaled $C^{1/2}$ norm bounded
by $Ct^{-1+\alpha/2}[f]_{C^\alpha}$. The constants in the Neumann Schauder estimate are controlled by the
common $C^3$ bounds for the boundary charts.
The estimate holds first for sufficiently small $t>0$.
It extends to the remaining times by the semigroup property.

By reflection in the common charts and a fixed partition of unity,
we obtain uniform extension bounds in $W^{2,p}$ and $C^{1,\gamma}$
for fixed $1<p<\infty$ and $0<\gamma<1$
\cite[Theorem~7.25 and Lemma~6.37]{GT}.
The Sobolev and Morrey inequalities and the interpolation estimates
hold with constants independent of $j$
\cite[Theorems~7.17 and~7.26 and Section~6.8]{GT}.
Let $u\in C^\infty(\bar\Omega_j)$ satisfy $u=\Phi_j$ on
$\partial\Omega_j$. For each fixed $1<p<\infty$,
\cite[Theorem~9.13]{GT} gives
\begin{equation}
 \|u\|_{W^{2,p}(\Omega_j)}\le C_p\{
 \|\Delta u\|_{L^p(\Omega_j)}+
 \|\Phi_j\|_{W^{2,p}(\Omega_j)}+
 \|u\|_{L^p(\Omega_j)}\}.
 \label{eq:uniform-dirichlet-w2p}
\end{equation}
Since $\Delta u=4\sum_i u_{i\bar i}$ and $\|\Phi_j\|_{C^3}\le C$,
\eqref{eq:uniform-dirichlet-w2p} gives
\[
 \|u\|_{W^{2,p}(\Omega_j)}
 \le C_p\bigl(1+\|\sum_i u_{i\bar i}\|_{L^p(\Omega_j)}
               +\|u\|_{L^p(\Omega_j)}\bigr).
\]
Uniform interior Schauder estimates also hold for the Poisson equation
and for linear uniformly elliptic equations on fixed rescaled balls
\cite[Section~6.1]{GT}.

\Needspace{5\baselineskip}
\begin{lemma}
\label{lem:weighted-trace-gradient}
Let $u$ be a smooth plurisubharmonic function on one of the domains
above, with the prescribed boundary value and a common $C^0$ bound.
For every finite $p>2n+1$,
\begin{equation}
 \|u\|_{C^{1,\delta_p}(\bar\Omega)}
 \le C_p\left[1+\left(\int_\Omega(-\rho)(\Delta u)^p\,dy\right)^{1/p}\right],
 \qquad \delta_p=1-\frac{2n+1}{p}.
 \label{eq:weighted-trace-gradient}
\end{equation}
The constants are uniform for the smooth approximations.
\end{lemma}
\begin{proof}
Let $w=u-\Phi$, where $\Phi$ is the prescribed boundary
extension. Then $w=0$ on $\partial\Omega$ and
$\Delta w=\Delta u-\Delta\Phi$.
Let $G_0$ be the Dirichlet Green function for $-\Delta$
\cite[Section~2.4]{GT}.
With $r=|x-y|$, symmetry and \cite[Theorem~2.3(i)]{Widman} give
$0\le G_0(x,y)\le Cr^{2-2n}\min\{1,d(y)/r\}$.
Applying the interior and boundary Schauder estimates
\cite[Sections~6.1--6.2]{GT} in $x$ on balls of radius comparable to $r$ gives
\begin{equation}
 |D_x^kG_0(x,y)|\le Cr^{2-2n-k}
                   \min\{1,d(y)/r\},\qquad k=1,2.
 \label{eq:laplace-green-derivatives}
\end{equation}%
These estimates use only the common $C^3$ geometry. Write $p'=p/(p-1)$.
On each annulus $A_r(x)=\Omega\cap\{r<|x-y|<2r\}$,
\eqref{eq:laplace-green-derivatives} gives
\begin{equation}
 \left(\int_{A_r(x)}|D_x^kG_0(x,y)|^{p'}
                      (-\rho(y))^{-1/(p-1)}\,dy\right)^{1/p'}
       \le C_pr^{2-k-(2n+1)/p}.
 \label{eq:weighted-green-annulus}
\end{equation}%
The comparison $-\rho\asymp d$ extends to all of $\Omega$ by
compactness, with constants uniform under smooth approximation.
On $d(y)\ge r$, we therefore have
$(-\rho(y))^{-1/(p-1)}\le C_pr^{-1/(p-1)}$.
On $d(y)<r$, the boundary factor leaves precisely
$d(y)^{p'-1/(p-1)}=d(y)$, whose annular integral is at most
$Cr^{2n+1}$.

For $x,x'\in\Omega$, set $\ell=|x-x'|$.
For $r\ge C\ell$ with $C$ fixed and sufficiently large, the points
$x,x'$ can be joined inside $\Omega\cap B_{r/2}(x)$ by a path of
length at most $C\ell$. Such paths follow from the fixed boundary
graph charts, and straight segments suffice in the interior.
Using \eqref{eq:laplace-green-derivatives} with $k=2$ and the same
calculation as in \eqref{eq:weighted-green-annulus}, we obtain
\begin{equation}
 \begin{aligned}
 &\left(\int_{A_r(x)}
 |D_xG_0(x,y)-D_xG_0(x',y)|^{p'}
 (-\rho(y))^{-1/(p-1)}\,dy\right)^{1/p'}\\
 &\qquad\le C_p\ell r^{-(2n+1)/p},\qquad r\ge C\ell.
 \end{aligned}
 \label{eq:weighted-green-difference}
\end{equation}
By Green's representation formula,
$Dw(x)=-\int_\Omega D_xG_0(x,y)\Delta w(y)\,dy$.
On $B_{C\ell}(x)\cup B_{C\ell}(x')$, we apply
\eqref{eq:weighted-green-annulus} with $k=1$ separately at $x$ and $x'$.
On the complement, we use \eqref{eq:weighted-green-difference}.
H\"older's inequality \cite[Section~7.1]{GT} and summation over dyadic
annuli give
\[
 \begin{aligned}
 |Dw(x)-Dw(x')|
 &\le C_p\left(\int_\Omega(-\rho)|\Delta w|^p\,dy\right)^{1/p}
 \left[\sum_{r\le C\ell}r^{1-(2n+1)/p}
       +\ell\sum_{r\ge C\ell}r^{-(2n+1)/p}\right]\\
 &\le C_p\left[1+\left(\int_\Omega(-\rho)(\Delta u)^p\,dy\right)^{1/p}\right]\ell^{\delta_p},
 \end{aligned}
\]%
since $\delta_p>0$.
The estimate for $D_xG_0$ also bounds $\|Dw\|_{L^\infty(\Omega)}$.
Restoring the boundary extension and using the common $C^0$ bound
completes the proof.
\end{proof}

We now check that the constants in Proposition~\ref{prop:elliptic-finish}
are uniform under smooth approximation.
By \eqref{eq:FEU1}, the operator in \eqref{eq:FEU3a} has common
ellipticity constants. The $C^{2,\alpha/2}$ norms of the boundary
graphs and boundary values are controlled by the common $C^3$ bounds. Using \cite[Theorem~1.4(b)]{SS}, we obtain
\eqref{eq:FEU4} with common $\gamma$ and $C_\gamma$.
Using the linear Schauder estimates \cite[Sections~6.1--6.2]{GT},
we fix $s>0$ independently of $j$ so that $Cs^\gamma\le1/2$ in
\eqref{eq:FEU9}. Absorbing the first term on the right gives
\eqref{eq:FEU2} with a uniform constant.

The following lemma passes the $C^{2,\alpha}$ estimate to the
limit on the original domain.

\begin{lemma}\label{lem:endpoint-limit}
Let $u_j$ be the smooth solutions of \eqref{intro:dirichlet} on $\Omega_j$ with
right-hand side $f_j$ and boundary value $\varphi_j$, as in
Lemma~\ref{lem:smooth-data}.
Suppose
\[
 0<\alpha<1,\qquad
 \sup_j\|u_j\|_{C^{2,\alpha}(\bar\Omega_j)}\le M.
\]
After extension to a fixed neighborhood of $\bar\Omega$, a subsequence
converges in $C^2(\bar\Omega)$ in the original complex coordinates.
The limit $u$ solves \eqref{intro:dirichlet} and satisfies
\[
 u\in C^{2,\alpha}(\bar\Omega),\qquad
 \|u\|_{C^{2,\alpha}(\bar\Omega)}\le CM,
\]
where $C$ depends only on the fixed geometry.
\end{lemma}

\begin{proof}
By the uniform extension estimates \cite[Lemma~6.37]{GT} and the
Arzel\`a--Ascoli theorem \cite[Theorem~7.25]{Rudin}, a subsequence of
the extensions converges in $C^2$ on a fixed compact neighborhood of
$\bar\Omega$ to $u$.
Every compact subset of $\Omega$ is contained in a connected
compact subset of $\Omega$ containing $z_0$. Uniform convergence of
$\rho_j$ makes $\rho_j<0$ on that set for all large $j$, so the set is
contained in the component $\Omega_j$. It follows that $(u_{i\bar j})\ge0$
and $\det(u_{i\bar j})=f$ in $\Omega$. For $p\in\partial\Omega$,
choose $p_j\in\partial\Omega_j$ with $p_j\to p$. Since
$u_j(p_j)=\Phi_j(p_j)$, we have $u(p)=\varphi(p)$.
For $x,y\in\Omega$, letting $j\to\infty$ gives
$|D^2u(x)-D^2u(y)|\le CM|x-y|^\alpha$.
By continuity this holds for $x,y\in\bar\Omega$.
\end{proof}

\section{Gradient and boundary estimates}\label{sec:gradient-boundary}

In this section, we prove global $C^{1,\beta}$ bounds for every fixed
$0<\beta<1$ and then estimate the full real Hessian on the boundary.

Unless otherwise stated, the domains are the smooth strongly pseudoconvex
approximations of Lemma~\ref{lem:smooth-data}. The boundary data are smooth,
and the geometric and boundary $C^3$ bounds are uniform. Let $u$ be a smooth
strictly plurisubharmonic solution of
\eqref{intro:dirichlet} with a smooth right-hand side $f$ satisfying
$0<\lambda\le f\le\Lambda$.

Chen--Cheng \cite[Remark~1.2]{CC} noted that their global gradient
estimate extends to the Dirichlet problem once global $C^0$ and
boundary gradient bounds are available. Lemma~\ref{lem:geometric-barriers}
provides these bounds under our $C^3$ assumptions on the strongly
pseudoconvex boundary and boundary values. We apply their result and
check that the constants remain uniform under smooth approximation.

\begin{proposition}\label{prop:uniform-gradient}
Let $0<\alpha<1$, $0<\lambda\le\Lambda$, and $F_*<\infty$.
Suppose that $u\in C^\infty(\bar\Omega)$ is a strictly
plurisubharmonic solution of \eqref{intro:dirichlet} on a
smooth bounded strongly pseudoconvex domain $\Omega$, with smooth
boundary value $\varphi$ and right-hand side satisfying
\[
 \lambda\le f\le\Lambda,\qquad \|f\|_{C^\alpha(\bar\Omega)}\le F_*.
\]
Then
\begin{equation}
 \|\partial u\|_{L^\infty}\le B.
 \label{eq:FE2}
\end{equation}
The constant $B$ depends only on $n$, $\alpha$, $\lambda$,
$\Lambda$, $F_*$, and the fixed $C^3$ geometric and boundary data
in Lemma~\ref{lem:geometric-barriers}. The same bound holds on the
smooth approximating domains of Lemma~\ref{lem:smooth-data} with
these common data.
\end{proposition}

\begin{proof}
By \eqref{eq:CG5} and the bounds for $f$, we have
\[
 \|u\|_{L^\infty}+\|Du\|_{L^\infty(\partial\Omega)}
       +\|\log f\|_{C^\alpha}\le C.
\]
The gradient estimate in \cite[Corollary~1.1 and Remark~1.2]{CC},
applied to the potential $u-|z|^2$ in the Euclidean metric, gives
\eqref{eq:FE2}.

We check uniformity on the approximating domains. The divergence
theorem gives
\begin{equation}
 \int_\Omega\sum_i u_{i\bar i}\,dy
 =\frac14\int_{\partial\Omega}u_\nu\,dS\le C.
 \label{eq:FE2b}
\end{equation}
The boundary gradient bound allows a common truncation level in the
energy argument. For each smooth solution, the truncated test
functions vanish near the boundary. The Euclidean $W^{1,1}$ Sobolev
inequality \cite[Theorem~7.10]{GT} applies to their zero extensions
with a dimensional constant. The additional threshold terms are
controlled by \eqref{eq:FE2b}. Together with the uniform extension
and smoothing bounds in Lemma~\ref{lem:smooth-data} and the fixed
geometric bounds, this makes the constants in the absorption and
iteration of \cite{CC} uniform.
\end{proof}

By tangential differentiation of $u=\Phi$ and \eqref{eq:FE2}, we obtain
bounds for the real tangential second derivatives. In a unitary frame
whose last vector is complex normal,
\eqref{eq:GC20} and the Schur complement formula
\cite[Section~0.8]{HJ} give
\begin{equation}
 \begin{gathered}
 (u_{i\bar j})=\begin{pmatrix}A_T&b_T\\b_T^*&c_N\end{pmatrix},
 \qquad c_TI\le A_T\le CI,\\
 c_N=\frac f{\det A_T}+b_T^*A_T^{-1}b_T
                         \le C(1+|b_T|^2).
 \label{eq:boundary-schur}
 \end{gathered}
\end{equation}
If the real mixed tangential--normal derivatives are bounded by $M$, then
$|b_T|\le C(1+M)$. Since $c_N>0$ and
$u_{\nu\nu}=4c_N-u_{J\nu,J\nu}$, the full real boundary Hessian is
bounded by $C(1+M^2)$.

\begin{lemma}\label{lem:cap-comparison}
Let $u$ be a smooth strictly plurisubharmonic solution with the fixed
positive determinant bounds and boundary value $\Phi$.
Let $\rho$ be the fixed strictly plurisubharmonic defining function,
negative inside.
There are fixed constants $K_0,D_0,\kappa>0$ and $0<R_0\le1$ with
the following property. Suppose $p\in\partial\Omega$, $0<R\le R_0$,
$0<E\le1$, $0<H\le\kappa E$, and
$|\gamma-\gamma(p)|\le CE$ on a set
$D\subset B_R(p)\cap\{0<-\rho<H\}$, with a fixed constant $C$.
Set
\[
 J=\Phi+(\gamma(p)+D_0E)\rho+(K_0+E/H)\rho^2,
 \qquad S=u-J.
\]
On $D$, we have
\begin{equation}
 \begin{gathered}
 cE(-\rho)\le S\le CE(-\rho),\qquad 0\le-S_\nu(p)\le CE,\\
 -L_uS\ge c\bigl(1+\sum_i u^{i\bar i}\bigr)+c(E/H)|\partial\rho|_u^2.
 \end{gathered}
 \label{eq:thin-mixed-comparison}
\end{equation}
The fixed positive definite summand of $(J_{i\bar j})$ can be chosen
to have an arbitrarily large fixed determinant. All constants depend only
on the fixed geometric and boundary data, their third-derivative
continuity moduli, the determinant bounds and the constant in $|\gamma-\gamma(p)|\le CE$,
independently of derivatives of the density.

In particular, for $0<e<1$ we may take $E=R^e$ and
$H=\kappa R^{2-e}$ provided the following estimate holds
throughout the cap.
\begin{equation}
 |\gamma(z)-\gamma(p)|\le CR^e
 \quad\text{for }z\in B_R(p),\quad 0<-\rho(z)<\kappa R^{2-e}.
 \label{eq:quotient-cap-modulus}
\end{equation}
The estimates also hold on anisotropic caps satisfying the same
quotient bound.
\end{lemma}

\begin{proof}
By Proposition~\ref{prop:tangential-block} and the Schur complement
formula, fixed $K_0$ and $R_0$ can be chosen so that
$\Phi+\gamma(p)\rho+K_0\rho^2$ is uniformly strictly plurisubharmonic on
$B_{2R_0}(p)$. The determinant lower bound can be made arbitrarily large by increasing
$K_0$ and decreasing $R_0$.
Since $0<-\rho<H$, we have
\[
 \begin{aligned}
 dd^cJ={}&dd^c(\Phi+\gamma(p)\rho+K_0\rho^2)
       +E(D_0+2\rho/H)dd^c\rho\\
       &+2(E/H)i\partial\rho\wedge\bar\partial\rho.
 \end{aligned}
\]
For $D_0>2$, the last two terms are nonnegative. Choosing the fixed
determinant lower bound sufficiently large and contracting with
$(u^{i\bar j})$ gives the differential inequality in
\eqref{eq:thin-mixed-comparison}.
Taking $D_0$ large relative to the constant in $|\gamma-\gamma(p)|\le CE$
and $\kappa$ sufficiently small, we also obtain
\[
 S=(-\rho)\{D_0E-(\gamma-\gamma(p))+K_0\rho+E\rho/H\}\asymp E(-\rho);
 \qquad -S_\nu(p)=D_0E\rho_\nu(p).
\]
\end{proof}

\Needspace{5\baselineskip}
The $C^{1,\beta}$ bound for $u$ gives $\gamma|_{\partial\Omega}\in C^\beta$.
The following lemma improves the exponent to $2\beta/(1+\beta)$
along curves in $\partial\Omega$ whose tangent vectors lie in the
complex tangent spaces.

\begin{lemma}
\label{lem:horizontal-slope-lift}
Let $u$ be PSH on the fixed strongly pseudoconvex domain, let
$u=\Phi$ on the boundary, and suppose
$\|u\|_{C^{1,\beta}(\bar\Omega)}\le B$ for $0<\beta<1$.
Put $c=2\beta/(1+\beta)$. Then
\begin{equation}
 |\gamma(p)-\gamma(q)|\le C\ell^c
 \label{eq:horizontal-slope-lift}
\end{equation}
whenever $p,q\in\partial\Omega$ can be joined by a $C^1$ curve of length
$\ell\le\ell_0$ whose tangent vectors $X$ satisfy $\partial\rho(X)=0$.
The constants depend only on the common $C^3$ geometry, $\Phi$, $B$ and
$\beta$.
\end{lemma}
\begin{proof}
Let $\omega(\ell)$ be the supremum of $|\gamma(p)-\gamma(q)|$ over
endpoints of such curves of length at most $\ell$.
The assumed gradient modulus gives
$\omega(\ell)\le B\ell^\beta$.
Fix $R=\eta\ell$, with $\eta>0$ fixed small. For all sufficiently
small $\ell$, depending on $\eta$, we have $\ell^2+\ell R\le CR^{2-c}$.
After a fixed quadratic holomorphic change of coordinates at the initial
point of the curve, we have
$\rho=\rho_{x_n}(0)\operatorname{Re}z_n+|z|^2+O(|z|^3)$.
We parametrize the curve by arc length and denote its velocity by $X$.
Since $\partial\rho(X)=0$, the curve satisfies $|z'|\le C\ell$
and $|z_n|\le C\ell^2$, and the complex-normal component of $X$
is $O(\ell)$. Let $\chi$ be a fixed nonnegative product mollifier supported in
$\{|\xi'|<1,\ |\xi_n|<1\}$, with $\int_{\mathbb C^n}\chi=1$.
We use scales $R$ in $z'$ and $R^{2-c}$ in $z_n$, with
\[
 \chi_R(z)=\frac{\chi(z'/R,z_n/R^{2-c})}{R^{2n+2-2c}}.
\]
Here $u_R=u*\chi_R$ and $\Phi_R=\Phi*\chi_R$ are the convolutions
with $\chi_R$. For the complex Hessian of $u_R$, we write
\[
 g^R=\bigl((u_R)_{i\bar j}(q',-C_0R^{2-c}+i\operatorname{Im}q_n)\bigr).
\]
For large fixed $C_0$, the kernel support at
$(q',-C_0R^{2-c}+i\operatorname{Im}q_n)$ lies in $\Omega$.
For each boundary projection $y$ of a support point, there is a $C^1$
curve $\zeta:[0,1]\to\partial\Omega$ satisfying
\[
 \zeta(0)=q,\qquad \partial\rho(\dot\zeta)=0,\qquad
 \operatorname{length}(\zeta)\le C_1R,\qquad |y-\zeta(1)|\le CR^{2-c}.
\]
To obtain $\zeta$, lift the straight segment from $q'$ to $y'$
to the boundary by solving $\partial\rho(\dot\zeta)=0$. In the normalized
coordinates this is an ordinary differential equation for
$\operatorname{Im}\zeta_n$, with velocity bounded by $C(\ell+R)|\dot\zeta'|$.
The lift therefore changes the imaginary coordinate by at most
$C(\ell R+R^2)$. The imaginary coordinate of $y$ differs from that of
$q$ by at most $CR^{2-c}$, since the support point has this property
and its distance to the boundary is $O(R^{2-c})$. The boundary graph
then gives the same bound for $|y-\zeta(1)|$.
By the $C^{1,\beta}$ bound and the definition of $\omega$, we have
on this support
\[
 |u-\Phi-\gamma(q)\rho|
       \le CR^{2-c}\{\omega(C_1R)+R^c\}.
\]%
The complex Hessian $g^R$ is positive semidefinite and satisfies
\[
 0\le g_{j\bar j}^{R}
      \le C\{1+\omega(C_1R)/R^c\},\quad j<n,
 \qquad 0\le g_{n\bar n}^{R}\le CR^{-2(1-c)}.
\]%
The Cauchy--Schwarz inequality gives
\[
 |g_{j\bar n}^{R}|
 \le CR^{c-1}\{1+\omega(C_1R)/R^c\}^{1/2}.
\]%
For $w=u-\Phi$, the identities
\[
 w_{x_jx_n}=4\operatorname{Re}w_{j\bar n}-w_{y_jy_n};\qquad
 w_{y_jx_n}=-4\operatorname{Im}w_{j\bar n}+w_{x_jy_n}
\]
give bounds for the real mixed derivatives in the normal and
complex tangential directions after convolution.
The correction is at most $CR^{c-1}$, since
$\rho_{y_n}=O(R^{2-c})$ and $|w_{y_n}|\le CR^c$
on the support.
Using the $C^{1,\beta}$ bound and the smoothing scale $R^{2-c}$
in the $z_n$ variable, we obtain
\[
 |\partial_a\partial_b(u_R-\Phi_R)|
 \le CR^{(2-c)(\beta-1)}=CR^{2c-2},
 \qquad a,b\in\{x_n,y_n\}.
\]%

Define the boundary approximation
\[
 b_R(q)=\frac{\partial_{x_n}(u_R-\Phi_R)(q',-C_0R^{2-c}+i\operatorname{Im}q_n)}
                  {\rho_{x_n}(q)}.
\]%
Then
\[
 |b_R(q)-\gamma(q)|\le C\{\omega(C_1R)+R^c+\ell^2\},
\]
and differentiation along the curve gives
\[
 |b_R'|\le C\left\{1+
 R^{c-1}\{1+\omega(C_1R)/R^c\}^{1/2}
                      +\ell R^{2c-2}\right\}.
\]
With $W(\ell)=\omega(\ell)/\ell^c$ and $\lambda=C_1\eta<1$,
we integrate along the curve and divide by $\ell^c$.
Since $c<1$ and $\ell^2R^{2c-2}/\ell^c=\eta^{2c-2}\ell^c$,
we obtain
\[
 W(\ell)\le C\lambda^cW(\lambda\ell)
              +C_\lambda\sqrt{1+C W(\lambda\ell)}+C_\lambda.
\]
The constant $C$ in the first term depends only on the original
$C^{1,\beta}$ bounds and the geometry. By Young's inequality, choosing
$\lambda$ sufficiently small gives
\[
 W(\ell)\le\tfrac12\lambda^{c-\beta}W(\lambda\ell)+C_\lambda.
\]
Using $W(\lambda^N\ell)\le B(\lambda^N\ell)^{\beta-c}$ and iterating,
we obtain
\[
 W(\ell)\le C_\lambda\sum_{j=0}^{N-1}
       (\tfrac12\lambda^{c-\beta})^j
       +2^{-N}B\ell^{\beta-c}.
\]
Letting $N\to\infty$ gives a uniform bound for $W$. The bound for
$\gamma$ covers the remaining curve lengths.
\end{proof}

\Needspace{5\baselineskip}
The next lemma estimates $(u-v)/\rho$ under a regularity assumption
on the reference solution $u$ alone.

\begin{lemma}
\label{lem:one-reference-quotient-comparison}
Let $u,v$ be smooth strictly PSH functions with the same fixed boundary
value and common positive determinant bounds. Suppose the
quotient $(u-\Phi)/\rho$ satisfies
\eqref{eq:quotient-cap-modulus} for one $0<e<1$, with the fixed data
used in Lemma~\ref{lem:cap-comparison}. Put
\[
 Q_0=\|u-v\|_{L^\infty},\qquad
 M_0=\|\det(v_{i\bar j})^{1/n}-\det(u_{i\bar j})^{1/n}\|_{L^\infty},
 \qquad \theta_e=\frac{2+(n-2)e}{2(n+1)-2e}.
\]
For $Q_0$ sufficiently small,
\begin{equation}
 \left\|\frac{u-v}{\rho}\right\|_{L^\infty}
 \le C M_0^{1-\theta_e}Q_0^{\theta_e}+CQ_0.
 \label{eq:one-reference-quotient-comparison}
\end{equation}
If instead only $u$ has a common global $C^{1,\beta}$ bound,
$0<\beta<1$, the same estimate holds without the hypothesis
\eqref{eq:quotient-cap-modulus}, with $\theta_e$ replaced by
\begin{equation}
 \vartheta_\beta^{\mathrm{q}}
 =\frac{1+(n-1)\beta}{n+1+(n-1)\beta},\qquad
 \vartheta_\beta^{\mathrm{q}}-\frac\beta{1+\beta}
 =\frac{1-\beta}{(1+\beta)\{n+1+(n-1)\beta\}}>0.
 \label{eq:holder-quotient-exponent}
\end{equation}
\end{lemma}
\begin{proof}
If $M_0=0$, uniqueness gives $u=v$.
Otherwise fix $p\in\partial\Omega$. On
$\Omega\cap B_R(p)\cap\{0<-\rho<\kappa R^{2-e}\}$,
Lemma~\ref{lem:cap-comparison} gives the comparison function
\[
 J=\Phi+(\gamma(p)+D_0R^e)\rho
       +(K_0+\kappa^{-1}R^{2e-2})\rho^2,\qquad
 S=u-J\asymp R^e(-\rho).
\]%
The fixed positive definite summand in the complex Hessian of $J$ has a large determinant.
The additional rank-one term in the complex normal direction has
coefficient at least $cR^{2e-2}$.
Thus, for $B/A$ below a fixed constant,
\begin{equation}
 \det\bigl(J_{i\bar j}-(B/A)I\bigr)^{1/n}
 \ge\det(u_{i\bar j})^{1/n}+cR^{-2(1-e)/n}.
 \label{eq:reference-determinant-gap}
\end{equation}%
Indeed, half the fixed positive matrix gives a determinant larger
than the prescribed upper density bound. Applying the matrix determinant lemma to the remaining positive
matrix and the rank-one normal term gives
\eqref{eq:reference-determinant-gap}.
Choose
\[
 B=C_1Q_0/R^2,\qquad
 A=C_2\{B+M_0R^{2(1-e)/n}\},
\]%
with fixed large $C_1,C_2$, and first assume $A\le1$.
The lower comparison function
\[
 W_-=(1-A)u+AJ-B|z-p|^2
\]
is PSH and has determinant at least that of $v$, by Minkowski's
determinant inequality \cite[Section~7.8]{HJ} and
\eqref{eq:reference-determinant-gap}.
We have $W_-\le v$ on the cap boundary: the quadratic term gives this
inequality on $|z-p|=R$, while $AS\ge Q_0$ gives it on $-\rho=\kappa R^{2-e}$,
and $u=v$ on $\partial\Omega$.
The comparison principle \cite[Theorem~A]{BT} therefore gives $W_-\le v$ throughout the cap.

For the upper comparison put $W_+=(1+A)u-AJ+B|z-p|^2$.
It is above $v$ on the cap boundary. If $v-W_+$ had a positive
interior maximum, then $(W_+)_{i\bar j}\ge v_{i\bar j}>0$ there.
At that point
\[
 (1+A)(u_{i\bar j})=(W_+)_{i\bar j}
                  +A\{J_{i\bar j}-(B/A)I\}.
\]
Minkowski's determinant inequality gives
$\det((W_+)_{i\bar j})^{1/n}
\le\det(u_{i\bar j})^{1/n}-cAR^{-2(1-e)/n}
<\det(v_{i\bar j})^{1/n}$, a contradiction.
This argument requires positivity of $W_+$ only at the contact point.
Thus $|u-v|\le AS+B|z-p|^2$ in the cap.

Along the inward normal segment from $p$, division by $-\rho$ gives
\[
 |u-v|/(-\rho)\le C\{AR^e+BR^{2-e}\}
       \le C\{M_0R^{e+2(1-e)/n}+Q_0/R^{2-e}\}.
\]%
For $-\rho\ge\kappa R^{2-e}$, the bound $|u-v|/(-\rho)\le CQ_0/R^{2-e}$ applies.

Choose $R=\min\{R_0,(Q_0/M_0)^{n/[2(n+1-e)]}\}$.
Then
\[
 A\le
 \begin{cases}
 C M_0^{n/(n+1-e)}Q_0^{(1-e)/(n+1-e)},&R<R_0,\\
 CQ_0,&R=R_0.
 \end{cases}
\]
Since $M_0$ is uniformly bounded, $A\le1$ for sufficiently small $Q_0$. The resulting estimate is
\eqref{eq:one-reference-quotient-comparison}.
For the assertion under a $C^{1,\beta}$ bound, set
$e=2\beta/(1+\beta)$ and use the same comparison function $J$ on
$\{0<-\rho<\kappa R^{2-e},\ |z'|<R,\ |\operatorname{Im}z_n|<\kappa R^{2-e}\}$.
Lift the straight segment from $0$ to $z'$ as in
Lemma~\ref{lem:horizontal-slope-lift}. Its endpoint has imaginary
normal coordinate $O(R^2)$ and lies within $O(R^{2-e})$ of the cap
point, by the boundary graph and the cap inequalities.
By Lemma~\ref{lem:horizontal-slope-lift} and the assumed
$C^{1,\beta}$ bound for $u$, we have
$|(u-\Phi)/\rho-\gamma(p)|\le C\{R^e+R^{\beta(2-e)}\}\le CR^e$.
Set
\[
 \begin{gathered}
 \chi=|z'|^2+\kappa^{-2}R^{2e-2}(\operatorname{Im}z_n)^2;\\
 B=CQ_0/R^2\quad\text{and}\quad A=C_\kappa[B+M_0R^{2(1-e)/n}].
 \end{gathered}
\]
Since $R\le1$, we have
\[
 dd^c\chi\le C_\kappa\{I+R^{2e-2}
 i\partial\rho\wedge\bar\partial\rho\}.
\]
Thus $dd^cJ-(B/A)dd^c\chi$ stays positive with the same
determinant gap, after increasing $C_\kappa$.
The preceding comparison arguments based on Minkowski's determinant inequality
therefore apply without change.
For a point on the inward normal segment from $p$ in the original coordinates,
with $-\rho<c\kappa R^{2-e}$, the normalized
coordinates satisfy $z',\operatorname{Im}z_n=O(-\rho)$, so
\[
 \chi/(-\rho)\le C_\kappa R^e,\qquad
 B\chi/(-\rho)\le C_\kappa Q_0/R^{2-e}.
\]
For $-\rho\ge c\kappa R^{2-e}$, the bound
$|u-v|/(-\rho)\le C_\kappa Q_0/R^{2-e}$ applies directly.
Thus \eqref{eq:one-reference-quotient-comparison} holds with
$e=2\beta/(1+\beta)$ and $\theta_e=\vartheta_\beta^{\mathrm{q}}$.
\end{proof}

Let $u$ be smooth and strictly PSH on the closure of a fixed smooth bounded
domain $\Omega$. Suppose $0<\lambda\le f=\det(u_{i\bar j})\le\Lambda$ and
$\|\partial u\|_{L^\infty}\le B$. Let $G_x\ge0$ be the positive
Dirichlet Green function \cite[Theorem~1.1]{GW}, with $G_x=0$ on $\partial\Omega$.
All derivatives of $G_x(y)$ below are taken in $y$. The normalization is
\begin{equation}
                 -f L_uG_x=\delta_x.
 \label{eq:FG1}
\end{equation}
By \eqref{eq:cofactor-divergence},
$fL_uh=\sum C^{i\bar j}h_{i\bar j}$. With $C=(C^{i\bar j})$
and real coordinates ordered as $(x,y)$, let
\begin{equation}
 A=\frac14\begin{pmatrix}\operatorname{Re}C&-\operatorname{Im}C\\ \operatorname{Im}C&\operatorname{Re}C\end{pmatrix},
 \qquad fL_uh=\operatorname{div}(A\nabla h).
 \label{eq:real-cofactor}
\end{equation}
Here $\operatorname{Re}C$ is symmetric and $\operatorname{Im}C$ is skew symmetric.
By \eqref{eq:cofactor-divergence}, every row of $A$ has zero divergence.
The second identity in \eqref{eq:real-cofactor} follows by expanding
$\partial_i\partial_{\bar j}$. Thus \eqref{eq:FG1} is the equation for the Green function of a
real symmetric operator in divergence form.
Let
$\lambda_i$ be the eigenvalues of $(u_{i\bar j})$.
The real eigenvalues of $A$ are $f/(4\lambda_i)$, each with
multiplicity two. Hence
\begin{equation}
 \det_{\mathbb R}A=4^{-2n}f^{2n-2},\qquad
                   \operatorname{tr}A^{-1}=2\Delta u/f\le C\Delta u.
 \label{eq:FG2}
\end{equation}

For $b\in W^{1,1}_0(\Omega)$ and $h\in L^\infty(\Omega)$,
symmetry of $p_t$ and $P_t1=1$ give
\begin{equation}
 \int_\Omega b(h-P_th)\,dy
 =\frac12\iint_{\Omega\times\Omega}p_t(y,z)
             (b(y)-b(z))(h(y)-h(z))\,dy\,dz.
 \label{eq:green-bilinear}
\end{equation}
The zero extension of $b$ satisfies \cite[Chapter~4]{EG}
$\|b(\cdot+\xi)-b\|_{L^1(\mathbb R^{2n})}
\le|\xi|\|Db\|_{L^1(\Omega)}$. Integrating this inequality against
the Gaussian bound in \eqref{eq:uniform-neumann-kernel}, we obtain
\[
 \left|\int_\Omega b(I-P_t)h\right|
 \le C\sqrt t\,\|Db\|_{L^1}\|h\|_{L^\infty},
 \qquad0<t\le1.
\]

By Subsection~\ref{sec:analytic-uniformity}, the Neumann heat semigroup
satisfies
\begin{equation}
 \|P_tf-f\|_{L^\infty}\le Ct^{\alpha/2}[f]_\alpha,
 \qquad\|\nabla P_tf\|_{L^\infty}
                       \le Ct^{(\alpha-1)/2}[f]_\alpha.
 \label{eq:neumann-heat-bounds}
\end{equation}
It follows from \eqref{eq:neumann-heat-bounds} that
\begin{equation}
 \sup_{0<t\le1}\|P_tf\|_{C^\alpha(\bar\Omega)}
                  \le C\|f\|_{C^\alpha(\bar\Omega)}.
 \label{eq:neumann-holder-bound}
\end{equation}
The Neumann heat regularization also satisfies
$\lambda\le P_tf\le\Lambda$ and $P_tf\in C^\infty(\bar\Omega)$ for every $t>0$.

\subsection{Auxiliary flows and global gradient estimates}

We will obtain global $C^{1,\beta}$ bounds for every fixed
$0<\beta<1$ from weighted integral estimates along an auxiliary
parabolic flow. Fix $0<\alpha<1$.

\begin{lemma}\label{lem:para-flow}
There are smooth auxiliary flows $\psi_\varepsilon$ satisfying the
corner compatibility conditions and the following common estimates
on a fixed short interval. Their initial values $u_\varepsilon$ satisfy
$\det((u_\varepsilon)_{i\bar j})=P_{\varepsilon^2}f$ and
$u_\varepsilon|_{\partial\Omega}=\varphi$. Put
\[
 g=((\psi_\varepsilon)_{i\bar j}),\quad Q=\det g,\quad U=\partial_t\psi_\varepsilon,
 \quad L=g^{i\bar j}\partial_i\partial_{\bar j},\quad
 \mathcal{T}=\operatorname{tr}g,\quad v=\operatorname{tr}g^{-1}.
\]
Here $\rho$ is a fixed strictly plurisubharmonic defining function.
The flows have boundary value $\varphi$, satisfy
$\partial_t\psi_\varepsilon=\log Q-H_\varepsilon$, and obey
\begin{equation}
 \begin{gathered}
 c\le Q\le C,\qquad \|D\psi_\varepsilon\|_{L^\infty}\le C,
 \qquad \int_\Omega \mathcal{T}\le C,\\
 |U(t,z)|\le Ct^{\alpha/2}\min\{1,-\rho(z)/t\},\qquad
 \|\psi_\varepsilon(t)-u_\varepsilon\|_{L^\infty}
       \le Ct^{1+\alpha/2},\\
 \|D^2H_\varepsilon(t)\|_{L^\infty}
       +\|H_{\varepsilon,t}(t)\|_{L^\infty}
       \le C[1+t^{-1+\alpha/2}].
 \end{gathered}\label{eq:para-flow-bounds}
\end{equation}
Moreover, $e^{H_\varepsilon(t)}$ has a uniform spatial $C^\alpha$ bound.
\end{lemma}

\begin{proof}
By \cite{CKNS}, the Dirichlet problem with density
$f_\varepsilon=P_{\varepsilon^2}f$ and boundary value $\varphi$ has a
smooth strictly plurisubharmonic solution $u_\varepsilon$. Choose a smooth cutoff $\chi$ equal to zero on
$[0,1/4]$ and one on $[1/2,\infty)$.
Define
\[
 H_\varepsilon(t)=\log\bigl[(1-\chi(t/\varepsilon^2))f_\varepsilon
       +\chi(t/\varepsilon^2)P_{t+\varepsilon^2}f\bigr].
\]
By \eqref{eq:uniform-neumann-kernel} and
\eqref{eq:neumann-holder-bound}, $H_\varepsilon$ has a common spatial
$C^\alpha$ bound and satisfies
\[
 \begin{gathered}
 \|H_{\varepsilon,t}\|_{L^\infty}+\|D^2H_\varepsilon\|_{L^\infty}
       \le C[1+(t+\varepsilon^2)^{-1+\alpha/2}];\\
 \|DH_\varepsilon\|_{L^\infty}
       \le C[1+(t+\varepsilon^2)^{(\alpha-1)/2}].
 \end{gathered}
\]
We have $H_\varepsilon(t)=\log f_\varepsilon$ for
$0\le t\le\varepsilon^2/4$.

For $A$ fixed sufficiently large, the time-independent strict
subsolution $\Phi+A\rho$ has determinant at least
$\sup_{t,z}e^{H_\varepsilon}$ and lies below $u_\varepsilon$ by
comparison. The smooth Dirichlet existence theorem
\cite[Theorem~1]{HouLi} therefore gives a smooth flow for each initial value
on any fixed finite interval. Since the forcing is independent of $t$
near zero and $u_\varepsilon$ solves the corresponding stationary
equation, all corner compatibility conditions hold.

By the parabolic maximum principle \cite[Chapter~II]{Lieberman}, applied
to $(\partial_t-L)U=-H_{\varepsilon,t}$ with zero initial and boundary
values, we obtain $|U|\le Ct^{\alpha/2}$.
It follows that $c\le Q=e^{U+H_\varepsilon}\le C$.
Integration in time also gives
$\|\psi_\varepsilon(t)-u_\varepsilon\|_{L^\infty}\le Ct^{1+\alpha/2}$. We increase
$A$ once more so that $\det((\Phi+A\rho)_{i\bar j})\ge\sup Q$.
At each time, comparison with this fixed barrier and the harmonic
extension of $\varphi$ gives the boundary gradient bound. Proposition~\ref{prop:uniform-gradient}
gives the initial gradient bound. Applying the maximum principle to
the real directional equations
$(\partial_t-L)\psi_e=-(H_\varepsilon)_e$ gives a common interior
gradient bound, since the displayed bound for $DH_\varepsilon$ is
integrable at zero. The divergence theorem then bounds $\int_\Omega \mathcal{T}$.

To estimate $U$ near the boundary, we compare each of $U$ and $-U$ with
$M t^{\alpha/2-1}\rho$ in the time-dependent region $0<-\rho<\kappa t$.
Since $L\rho\ge c\,v\ge c_0>0$, we choose $\kappa$ small and then
$M$ large. The comparison function has $(\partial_t-L)$-value at most
$-C t^{\alpha/2-1}$, is zero on the boundary, and is below
both $U$ and $-U$ on $-\rho=\kappa t$ by the preceding uniform bound.
The comparison region is empty at time zero, and $U$ tends to zero
at the initial corner.
Thus $|U|\le Ct^{\alpha/2-1}(-\rho)$ on $0<-\rho<\kappa t$.
Together with $|U|\le Ct^{\alpha/2}$, this gives the bound for $U$
in \eqref{eq:para-flow-bounds}.
\end{proof}

\Needspace{5\baselineskip}
The next corollary gives the rate of convergence of $(\psi_\varepsilon(t)-u_\varepsilon)/\rho$ as $t\downarrow0$.

\begin{corollary}
\label{cor:flow-quotient-approximation}
In any complex dimension $n\ge2$, suppose the initial functions
$u_\varepsilon$ in Lemma~\ref{lem:para-flow} satisfy
\eqref{eq:quotient-cap-modulus} uniformly with exponent $0<e<1$. Then
\begin{equation}
 \left\|\frac{\psi_\varepsilon(t)-u_\varepsilon}{\rho}\right\|_{L^\infty}
       \le C t^{\alpha/2+\theta_e}.
 \label{eq:flow-quotient-approximation}
\end{equation}
Here $\theta_e$ is given in Lemma~\ref{lem:one-reference-quotient-comparison}.
If the initial functions instead have a common $C^{1,\beta}$ bound,
the estimate holds with $\theta_e$ replaced by
$\vartheta_\beta^{\mathrm{q}}$ from \eqref{eq:holder-quotient-exponent}.
Neither assertion requires a spatial H\"older bound for the
determinant of the flow at positive times.
\end{corollary}
\begin{proof}
By the semigroup identity $P_{t+\varepsilon^2}f=P_tf_\varepsilon$,
the definition of $H_\varepsilon$ and the uniform $C^\alpha$ bound give
$\|H_\varepsilon(t)-\log f_\varepsilon\|_{L^\infty}\le Ct^{\alpha/2}$.
Together with \eqref{eq:para-flow-bounds} and $Q=e^{H_\varepsilon+U}$,
this yields
\[
 \|Q(t)^{1/n}-f_\varepsilon^{1/n}\|_{L^\infty}\le Ct^{\alpha/2}.
\]
Using $\|\psi_\varepsilon(t)-u_\varepsilon\|_{L^\infty}\le Ct^{1+\alpha/2}$
from \eqref{eq:para-flow-bounds},
Lemma~\ref{lem:one-reference-quotient-comparison} gives both estimates,
since $\theta_e<1$ and $\vartheta_\beta^{\mathrm{q}}<1$. The common
gradient bounds cover the remaining time interval.
\end{proof}

\begin{lemma}\label{lem:para-boundary-trace}
Fix $n\ge2$, an integer $p\ge2$, and $0<2\eta<\alpha$.
For every $\epsilon>0$, the flows of Lemma~\ref{lem:para-flow}
satisfy, at each smooth time slice,
\begin{equation}
 \int_{\partial\Omega}\mathcal{T}^p\,dS
 \le\epsilon\int_\Omega Q((L\rho)+(-\rho)^{\eta-1}a_N)\mathcal{T}^p+C_{p,\epsilon},
 \qquad a_N=|\partial\rho|_g^2.
 \label{eq:para-boundary-trace}
\end{equation}
The constants are independent of spatial H\"older norms of $Q(t)$.
\end{lemma}

\begin{proof}
At a boundary point, holomorphic linear and quadratic changes of
coordinates give $\rho=\rho_{x_n}(0)\operatorname{Re}z_n+|z|^2+O(|z|^3)$.
Let $K_j$, $1\le j\le2n-1$, be the infinitesimal generators of affine
unitary transformations preserving $\rho_{x_n}(0)\operatorname{Re}z_n+|z|^2$,
and let $N=\partial_{x_n}$.
On a chart of radius $r$, put
\[
 a_j=\frac{K_j\rho}{N\rho},\qquad T_j=K_j-a_jN.
\]
These fields span the real tangent spaces to the level sets of the
defining function and satisfy
$T_j\rho=0$, $a_j=O(r^2)$, $Da_j=O(r)$ and $D^2a_j=O(1)$.
The coordinate changes have common bounds. We include the
logarithm of the real Jacobian determinant in $H_\varepsilon$.

At a fixed time, let $\gamma=(\psi_\varepsilon-\Phi)/\rho$ and set
\[
 w_j=T_j(\psi_\varepsilon-\Phi)=\rho\alpha_j,\qquad
 \alpha_j=T_j\gamma.
\]
The common gradient bound implies $|\gamma|+|w_j|\le C$.
By direct calculation, we have
\begin{equation}
 -Lw_j=-T_j\log Q+R_j,\qquad
 \bigl(\sum_j|R_j|^2\bigr)^{1/2}\le Cv+Cr\sqrt v\Big(\sum_j|\partial w_j|_g^2\Big)^{1/2}.
 \label{eq:para-nearhol-source}
\end{equation}
To control $\partial(N\psi_\varepsilon)$, put
$a_c=JN\rho/N\rho$. Then
\[
 (i-a_c)\partial_i(N\psi_\varepsilon)
 =\partial_i\bigl((JN-a_cN)\psi_\varepsilon\bigr)
       +(N\psi_\varepsilon)\partial_i a_c+2i g_{i\bar n}.
\]
We express $JN-a_cN$ in the $T_j$ frame and take complex conjugates.
Contraction with $g^{i\bar k}(a_j)_i$ then reduces the term containing
$g_{n\bar k}$ to $(a_j)_n$. The remaining derivative has coefficient $Da_j=O(r)$,
giving \eqref{eq:para-nearhol-source}. The term $L(T_j\Phi)$ involves at most three derivatives of
$\Phi$ and $\rho$.

Using $d^c\rho$, we define
\[
 \begin{gathered}
 t_j=\frac{d^c\rho(T_j)}{\big(\sum_k|d^c\rho(T_k)|^2\big)^{1/2}},
 \qquad T_c=\sum_k t_kT_k,\qquad \tau=\sum_jt_j\alpha_j,\\
 h_j=\alpha_j-\tau t_j,\qquad X_j=T_j-t_jT_c,
 \qquad w_c=\sum_jt_jw_j=\rho\tau.
 \end{gathered}
\]
Then $\sum_jt_j^2=1$, $\sum_jh_jt_j=0$ and $X_j\gamma=h_j$.
The fields $X_j$ are complex tangential, and the coefficients $t_j$
have common $C^2$ bounds.
Hermitian pairings are summed componentwise for vector-valued functions.
Let $A$ be the real symmetric matrix in \eqref{eq:real-cofactor}, with
$C^{i\bar j}=Qg^{i\bar j}$, and set $\mathcal N_A h=A\nabla h\cdot\nu$.
Put $\mathcal B=\sum_jh_j\partial t_j$.
Direct differentiation, without assuming that $T_c$ is real holomorphic,
gives
\begin{align}
 -L(-\rho|h|^2)
 & =-(L\rho)|h|^2+2\rho\sum_j|\partial h_j|_g^2
       -4\operatorname{Re}\langle\partial w_c,\mathcal B\rangle_g\notag\\
 &\quad-2\rho\tau\sum_jh_jLt_j+2\sum_jh_j(T_j\log Q-R_j),
       \label{eq:para-horizontal-product}\\
 -Lw_c
 &=-T_c\log Q+\sum_jt_jR_j-\rho\sum_jh_jLt_j
       -\rho\tau\sum_j|\partial t_j|_g^2\notag\\
 &\quad-2\rho\operatorname{Re}\sum_j\langle\partial t_j,\partial h_j\rangle_g
       -2\operatorname{Re}\langle\partial\rho,\mathcal B\rangle_g.
       \label{eq:para-characteristic-product}
\end{align}%
In particular
\[
 \partial w_j=\rho(\partial h_j+\tau\partial t_j)
                   +h_j\partial\rho+t_j\partial w_c.
\]
Using $a_N\le Cv$, we obtain
\[
 |\partial w|_g\le C\{(-\rho)|\partial h|_g+
 \sqrt v|w_c|+\sqrt{a_N}|h|+|\partial w_c|_g\}.
\]

Fix an integer $p\ge2$, choose one sufficiently large $\kappa$, and
put $Z=|h|^2+\kappa (-\rho)\tau^2$. Then
\begin{equation}
 \begin{aligned}
 -L((-\rho)Z^p)
 ={}&-pZ^{p-1}L((-\rho)Z)-(p-1)(L\rho)Z^p\\
 &-p(p-1)(-\rho)Z^{p-2}|\partial Z|_g^2.
 \end{aligned}
 \label{eq:para-horizontal-power}
\end{equation}
After absorbing the geometric error terms, we obtain
\begin{equation}
 -L((-\rho)Z)\le
 -c\{v|h|^2+(-\rho)|\partial h|_g^2+\kappa|\partial w_c|_g^2\}
 +Cv+2\sum_jh_jT_j\log Q-2\kappa w_cT_c\log Q.
 \label{eq:para-geometric-base}
\end{equation}
Indeed, $|\mathcal B|_g^2\le Cv|h|^2$ controls the cross term for large
$\kappa$, while the commutator leaves $Crv|h|^2$, absorbed for small
$r$. The other frame terms use the bounded $\rho\tau=w_c$ and Young's
inequality. The remaining quadratic term is
$-2\kappa w_c \rho\tau\sum_j|\partial t_j|_g^2
=-2\kappa w_c^2\sum_j|\partial t_j|_g^2\le0$.
For $p\ge2$, the term $-(p-1)(L\rho)Z^p$ in
\eqref{eq:para-horizontal-power} controls the full quantity $vZ^p$.
Fix $0<2\eta<\alpha$. We use the weight
\[
 \mu=-\rho+(-\rho)^{1+\eta},\qquad
 a=\mu/(-\rho)=1+(-\rho)^\eta,
 \qquad L\rho\asymp v.
\]
The product rule gives
\begin{equation}
 \begin{aligned}
 -L(\mu Z^p)
 ={}&-aL((-\rho)Z^p)+\eta(-\rho)^\eta (L\rho)Z^p\\
 &-\eta(1+\eta)(-\rho)^{\eta-1}a_NZ^p\\
 &+2p \eta(-\rho)^\eta Z^{p-1}\operatorname{Re}
                \langle\partial\rho,\partial Z\rangle_g,
 \end{aligned}
 \label{eq:para-convex-distance}
\end{equation}
By Young's inequality, the last term is bounded by a small multiple
 of $(-\rho)Z^{p-2}|\partial Z|_g^2$ and
 $C_p(-\rho)^{2\eta-1}a_NZ^p$. Since $(-\rho)^\eta\to0$
 at the boundary, reducing the collar width absorbs this remainder
 and the positive term $\eta(-\rho)^\eta (L\rho)Z^p$.
 The retained terms include $(L\rho)Z^p$ and $(-\rho)^{\eta-1}a_NZ^p$.

We use a fixed finite cover by boundary charts and cutoffs $\chi_l$
such that the sets $\{\chi_l=1\}$ cover a fixed collar of the boundary.
Let $d_0>0$ be small and choose a smooth cutoff $0\le\sigma\le1$
equal to one on $[0,d_0/2]$ and zero on $[d_0,\infty)$.
Set $\theta_l=\sigma(-\rho)\chi_l$ and $\omega_{l,p}=\theta_l^{2p}$.
We write $Z_l,h_l$ for the quantities in the $l$th chart.
All fields tangent to the level sets annihilate $a$.
Write $A_T$ for the complex tangential block of $g$.
The elementary matrix inequalities needed below are
\begin{equation}
 a_N\mathcal{T}\ge c>0,\qquad (L\rho)\ge c(1+\mathcal{T}^{1/(n-1)}),\qquad
 |X_jf|^2\le C\operatorname{tr}A_T|\partial f|_g^2,
 \quad |T_cf|^2\le C\mathcal{T}|\partial f|_g^2.
 \label{eq:para-joint-matrix}
\end{equation}
The first inequality follows from the Cauchy--Schwarz inequality.
For the second, we use the arithmetic--geometric mean inequality
and the upper bound for $\det g$. Using the mixed Hessian integral estimate
with $j=n-1$ in the proof of Lemma~\ref{lem:boundary-quotient-base},
we also obtain $\int Q(L\rho)\le C$, since $\|\rho\|_{C^2}\le C$.

The flow estimate gives $|U|+|Q-e^{H_\varepsilon}|\le C(-\rho)^{\alpha/2}$.
Convolving $e^{H_\varepsilon}$ tangentially at radius $\sqrt{-\rho}$ and
adding $C_0(-\rho)^{\alpha/2}$ defines $F_g$. For a fixed sufficiently large
$C_0$, the difference $b=F_g-Q$ satisfies
\[
 0\le b\le C(-\rho)^{\alpha/2},\qquad |X_jF_g|\le C(-\rho)^{(\alpha-1)/2}.
\]
The added term is constant on each level set.

Let $\mathcal H=\sum_jh_jX_j$. The quadratic form
$\sum_jX_j\otimes X_j$ and its $J$-invariant part are uniformly
positive on the complex tangent spaces to the level sets of $\rho$. The identity
\[
 (-\rho)\gamma_{a\bar b}=\Phi_{a\bar b}+\gamma\rho_{a\bar b}-(A_T)_{a\bar b}
\]
therefore bounds the contraction of the invariant part with $D^2\gamma$ by
$-c\operatorname{tr}A_T/(-\rho)+C/(-\rho)$. The complementary part of the quadratic
form is expressed in the $X_j$ frame with bounded coefficients. Using
$X_kh_j=D^2\gamma(X_k,X_j)+(D_{X_k}X_j)\cdot D\gamma$,
\eqref{eq:para-joint-matrix} and
$\rho D\gamma=D(\psi-\Phi)-\gamma D\rho$, we obtain
\begin{equation}
 \operatorname{div}\mathcal H
 \le-c\operatorname{tr}A_T/(-\rho)+C/(-\rho)
                 +C\sqrt{\operatorname{tr}A_T}|\partial h|_g.
 \label{eq:para-source-divergence}
\end{equation}
The terms arising from derivatives of the frame are controlled
using $|D\gamma|\le C/(-\rho)$.

Since $Q\sum_jh_jT_j\log Q=\mathcal H Q$, integration by parts gives
\[
 -2p\int\omega_{l,p}aZ_l^{p-1}\mathcal H b
   =2p\int b\operatorname{div}
                     (\omega_{l,p}aZ_l^{p-1}\mathcal H).
\]
There is no boundary term, and $X_ja=0$.
By \eqref{eq:para-source-divergence}, the terms involving $b$ satisfy
\[
 \int\omega_{l,p}(-c_pb+C_pb^2)
       \frac{\operatorname{tr}A_T}{(-\rho)}Z_l^{p-1}\le0
\]
after reducing the collar. We retain the corresponding squared gradient terms. Derivatives of the cutoff give
$C_pb\theta_l^{2p-1}Z_l^{p-1/2}$, controlled by the retained
zeroth-order term and $C_p\int Q(L\rho)$.
Writing $Y_l=\theta_l^2Z_l$, the only singular remainder is
\[
 C_p\sum_l\int Q(-\rho)^{\alpha/2-1}Y_l^{p-1}.
\]
Using $(L\rho)|h|^2Z^{p-1}$, we bound the contribution from $F_g$ by
\[
 C_p(-\rho)^{\alpha-1}(L\rho)^{-1}Y_l^{p-1}
 \le C_p(-\rho)^{\alpha/2-1}Y_l^{p-1}.
\]%
Since $1/a_N\le C\mathcal{T}$, Young's inequality gives, for every $\delta>0$,
\begin{equation}
 (-\rho)^{\alpha/2-1}Y^{p-1}
 \le\delta (-\rho)^{\eta-1}a_NY^p
  +C_{p,\delta}(-\rho)^{p(\alpha/2-\eta)}(-\rho)^{\eta-1}a_N\mathcal{T}^p.
 \label{eq:para-joint-normal}
\end{equation}
The first term is absorbed for small $\delta$. Since
$p(\alpha/2-\eta)>0$, reducing the collar makes the coefficient of
the second arbitrarily small.

For the term containing $T_cQ$, we integrate by parts in $T_c$.
Using $|w_c|\le C$, $w_c^2/(-\rho)\le Z/\kappa$ and
\eqref{eq:para-joint-matrix}, we obtain small multiples of the retained
gradient terms and
$C_p\int Q\mathcal{T}Y_l^{p-1}$. For every $\delta>0$,
\begin{equation}
 \mathcal{T}Y^{p-1}\le\delta (L\rho)(Y^p+\mathcal{T}^p)+C_{p,\delta}(L\rho).
 \label{eq:para-joint-reaction}
\end{equation}
Indeed, division by $(L\rho)$ and Young's inequality leave $(\mathcal{T}/(L\rho))^p$.
Since $(L\rho)\ge c\mathcal{T}^{1/(n-1)}$, we have
\[
 (\mathcal{T}/(L\rho))^p\le\delta \mathcal{T}^p+C_{p,\delta}.
\]%
Bounded $w_c$ controls the frame and cutoff terms. Also $T_c\rho=0$
eliminates derivatives of $a$.

The diffusion term contributes $\mu Z_l^pQL\omega_{l,p}$ when the
operator acts on the cutoff.
Where tangential cutoffs vary, it is bounded by
$C_pd_0Q(L\rho)\chi_l^{2p-2}Z_l^p$. On overlaps, uniformly bounded changes of frame give
$|h_l|\asymp|h_m|$ and $|\tau_l|\le C(|\tau_m|+|h_m|)$, hence
$Z_l\asymp Z_m$. Since the sets $\{\chi_l=1\}$ cover the collar,
\[
 \sum_l\chi_l^{2p-2}Z_l^p\le C_p\sum_l\chi_l^{2p}Z_l^p.
\]
The retained $(L\rho)Z_l^p$ terms therefore absorb these errors for small
$d_0$. Where the normal cutoff varies, $-\rho$ is bounded below, so
$Z_l\le C(d_0)$ and $\int Q(L\rho)\le C$ bounds the remaining errors.

Since $\mu=0$ and $D\mu=-D\rho$ on $\partial\Omega$, the boundary flux is
\[
 B_p=\sum_l\int_{\partial\Omega}\chi_l^{2p}
       (\mathcal N_A\rho)|h_l|^{2p}\,dS.
\]
Choosing the coefficients in \eqref{eq:para-joint-normal} and
\eqref{eq:para-joint-reaction}, then reducing the collar, absorbs
the errors into the retained terms. Discarding the remaining
nonnegative terms gives
\[
 B_p\le\epsilon\int_\Omega Q((L\rho)+(-\rho)^{\eta-1}a_N)\mathcal{T}^p+C_{p,\epsilon}.
\]
The integrations are justified on $-\rho>\delta$ and then by
$\delta\downarrow0$: $\mu$ is $C^1$ on the closure and its Hessian
is locally integrable.

By Proposition~\ref{prop:tangential-block}, the complex tangential
block of $g$ is uniformly positive and bounded on $\partial\Omega$.
On $\partial\Omega$,
\[
 D^2(\psi_\varepsilon-\Phi)(X_j,\nu)
       =\rho_\nu h_j+\gamma X_j\rho_\nu.
\]
Thus the mixed block is bounded by $C(1+|h|)$, and the Schur
complement formula gives $\mathcal T\le C(1+|h|^2)$.
In a unitary frame adapted to the boundary,
$\mathcal N_A\rho=(\rho_\nu/4)\det A_T\asymp1$.
Hence \eqref{eq:para-boundary-trace} follows.
\end{proof}

\begin{lemma}\label{lem:para-trace-nash}
Let $0\le\beta<1$. Suppose the initial functions $u_\varepsilon$
have a common $C^{1,\beta}$ bound when $\beta>0$, and a common
Lipschitz bound when $\beta=0$. Suppose also that
\begin{equation}
 \begin{gathered}
 \|\psi_\varepsilon(t)-u_\varepsilon\|_{L^\infty}\le Ct^{q_0},
 \qquad q_0>1,\\
 \| (\psi_\varepsilon(t)-u_\varepsilon)/(-\rho)\|_{L^\infty}\le Ct^c,
 \qquad c\ge0.
 \end{gathered}
 \label{eq:para-initial-errors}
\end{equation}
Set $\gamma_\beta=(1-\beta)/(1+\beta)$ and
$m=\max\{\gamma_\beta,1-2c,0\}$.
For every finite integer $p\ge2$,
\begin{equation}
 Y_p(t)+\int_t^{2t}\!\int_\Omega Q(L\rho)\mathcal{T}^p
       \le C_p(1+t^{-(p-1)m}),
 \qquad Y_p(t)\asymp\int_\Omega (-\rho)\mathcal{T}^p.
 \label{eq:para-trace-rate}
\end{equation}
At $\beta=0$, the common boundary values and gradient bounds allow
$c=0$, hence $m=1$.
\end{lemma}

\begin{proof}
We adapt the weighted energy estimate of
Li--Wang--Zheng \cite[Lemma~3.13]{LWZ} by introducing a spatial weight
comparable to the boundary distance and retaining the boundary flux.
For $0<2\eta<\alpha$, set
\[
 \mu=\int_0^{-\rho} e^{-s^\eta}\,ds,\qquad
 h(U)=e^{U^2-U/2},\qquad Y_p=\int_\Omega Qh(U)\mu \mathcal{T}^p.
\]
Then $\mu\asymp-\rho$ and
\[
 \begin{gathered}
 L\mu=-e^{-(-\rho)^\eta}\bigl(L\rho+\eta(-\rho)^{\eta-1}a_N\bigr);\\
 |\partial\mu|_g^2/\mu\le e^{-(-\rho)^\eta}a_N/(-\rho).
 \end{gathered}
\]
Commuting mixed derivatives and differentiating the determinant, we obtain
\[
 (\partial_t-L)\mathcal{T}^p\le
 C_p\|D^2H_\varepsilon\|_{L^\infty}\mathcal{T}^p
       -p^2\mathcal{T}^p|\partial\log \mathcal{T}|_g^2,
 \qquad Q_t=QLU.
\]
In the resulting energy identity, the nonnegative quadratic form
in the gradients is exactly
\[
 |p\partial\log \mathcal{T}+2U\partial U|_g^2
                         +\tfrac74|\partial U|_g^2.
\]
Let $B(U)=\int_0^U h(s)\,ds$. After integration by parts, the terms
containing derivatives of $\mu$ are
\[
 \int_\Omega Q\mathcal{T}^p(h+B)L\mu
 +\int_\Omega QB\operatorname{Re}
       \langle\partial\mu,\partial \mathcal{T}^p\rangle_g
 +\int_{\partial\Omega}\mathcal{T}^p\mathcal N_A\rho.
\]
Since $B=O(U)$ and $h+B>0$ for small $|U|$, Young's inequality
bounds the second term by a small multiple of the gradient quadratic form
and $C\int_\Omega QU^2\mathcal{T}^pe^{-(-\rho)^\eta}a_N/(-\rho)$.
Using \eqref{eq:para-flow-bounds}, we have $U^2\le Ct^{\alpha-\eta}(-\rho)^\eta$.
The negative term $-\eta(-\rho)^{\eta-1}e^{-(-\rho)^\eta}a_N$
therefore absorbs the
remainder on a common short interval. Since $\mu=U=0$ on
$\partial\Omega$, no conormal derivatives of $\mathcal{T}$ or $U$ remain.
Lemma~\ref{lem:para-boundary-trace} allows us to absorb the boundary
flux into $\int_\Omega Q((L\rho)+(-\rho)^{\eta-1}a_N)\mathcal{T}^p$, up to a fixed constant.
It follows that
\begin{equation}
 \begin{gathered}
 Y_p'+c_pD_p\le C_p(1+t^{-1+\alpha/2})Y_p+C_p;\\
 D_p=\int_\Omega Q\left[
 ((L\rho)+(-\rho)^{\eta-1}a_N)\mathcal{T}^p+(-\rho)\mathcal{T}^{p-2}|\partial \mathcal{T}|_g^2\right].
 \end{gathered}
 \label{eq:para-trace-energy}
\end{equation}

Write $J_p=\int_\Omega (-\rho)\mathcal{T}^p$. By mollification followed by a correction of the boundary values,
we obtain $u_R=\varphi$ on $\partial\Omega$ and
\[
 \|u_R-u_\varepsilon\|_{L^\infty}
 +R\|D(u_R-u_\varepsilon)\|_{L^\infty}
 +R^2\|D^2u_R\|_{L^\infty}\le CR^{1+\beta}.
\]
We convolve a controlled extension of $u_\varepsilon$ and subtract
the boundary error of the convolution, extended by the normal
projection with a collar cutoff. The common
$C^3$ geometry controls this correction also when $\beta=0$.
For $z=\psi_\varepsilon-u_R$, choose upper bounds
\[
 E\ge\|z\|_{L^\infty},\quad F\ge\|z/(-\rho)\|_{L^\infty},\quad
 M\ge\|D^2u_R\|_{L^\infty}.
\]
Set $Z_0=\int_\Omega(-\rho)\mathcal T^{p-1}|\partial z|_0^2$.
Using $z^2\le F^2(-\rho)^2$, $g\le\mathcal TI$, the bounds for $Q$
and integration by parts, we obtain
\[
 Z_0\le C_p[EJ_p+EMJ_{p-1}+F^2J_{p-1}+E^2D_p].
\]
To estimate $J_p$, we integrate the $\Delta z$ term by parts.
The term containing $D\rho$ is integrated by parts once more.
Using $z=0$ on $\partial\Omega$ and $|z|\le F(-\rho)$, we obtain the bound
$CFJ_{p-1}+C_pF\sqrt{D_pJ_{p-1}}$ for the $D\rho$ term.
The $D\mathcal{T}$ term is bounded by $C_p\sqrt{D_pZ_0}$.
Young's inequality then gives
\begin{equation}
 J_p\le C_p\{ED_p+(M+F+F^2/E)J_{p-1}\}.
 \label{eq:para-intrinsic-nash}
\end{equation}
All integrals involving lower powers of $\mathcal{T}$ in this estimate
carry the weight $(-\rho)$.

For $R=(\epsilon t)^{1/(1+\beta)}$, \eqref{eq:para-initial-errors}
and $q_0>1$ allow us to take $E=C\epsilon t$ on a shorter common
interval, with
\[
 M+F+F^2/E\le C_\epsilon(1+t^{-m}).
\]
Since $J_1\le C$, H\"older's inequality with the measure $(-\rho)\,dx$ gives
$J_{p-1}\le C J_p^{(p-2)/(p-1)}$. By Young's inequality in
\eqref{eq:para-intrinsic-nash},
\[
 J_p\le C_p\epsilon tD_p+C_{p,\epsilon}(1+t^{-(p-1)m}).
\]
Since $\epsilon>0$ is arbitrary and $t+t^{\alpha/2}\to0$,
\eqref{eq:para-trace-energy} and $Y_p\asymp J_p$ give, for any fixed
$A>0$ on a shorter common interval,
\[
 Y_p'+(A/t)Y_p\le C_{A,p}t^{-1-(p-1)m}+C_p.
\]
For $A>(p-1)m$, integration after multiplication by $t^A$ gives
$Y_p\le C_p(1+t^{-(p-1)m})$. The initial term vanishes because
$t^AY_p(t)\to0$ for each smooth initial solution. No common bound for $Y_p$ at time zero is used.
Integration of \eqref{eq:para-trace-energy} on $[t,2t]$ proves
\eqref{eq:para-trace-rate}.
\end{proof}

\begin{proposition}\label{prop:para-all-gradient}
In every complex dimension $n\ge2$, for every $0<\alpha<1$ and every fixed
$0<\beta<1$, the smooth Dirichlet solutions with the prescribed
positive $C^\alpha$ density bounds satisfy
\begin{equation}
                  \|u\|_{C^{1,\beta}(\bar\Omega)}\le C_\beta.
 \label{eq:para-all-gradient}
\end{equation}
The constants are uniform for the controlled $C^3$ approximations.
\end{proposition}

\begin{proof}
Suppose a common stationary
$C^{1,\beta}$ bound is known for one $0\le\beta<1$, with $\beta=0$
interpreted as the common Lipschitz bound. For $\beta>0$,
Corollary~\ref{cor:flow-quotient-approximation} gives the quotient
error exponent $c=\alpha/2+\vartheta_\beta^{\mathrm{q}}
>\beta/(1+\beta)$. At $\beta=0$ we take $c=0$. In both cases,
Lemma~\ref{lem:para-trace-nash} applies with $m=\gamma_\beta$.
Thus Lemma~\ref{lem:weighted-trace-gradient} gives, for every finite
integer $p>2n+1$,
\[
 \|\psi_\varepsilon(t)\|_{C^{1,\delta_p}}
 \le C_pt^{-\gamma_\beta(1-1/p)}.
\]
For $t_j=2^{-j}t_0$, we have
\[
 \|\psi_\varepsilon(t_j)-\psi_\varepsilon(t_{j+1})\|_{L^\infty}
       \le Ct_j^{1+\alpha/2},
\]
with $C$ independent of $\varepsilon$. By H\"older interpolation
\cite[Lemma~3.2, equation~(3.10)]{LWZ}, the series
of these differences converges in $C^{1,\beta'}$ provided
\[
 (1+\alpha/2)(\delta_p-\beta')
    >\gamma_\beta(1-1/p)(1+\beta').
\]%
For any $\beta'$ satisfying
\[
 \beta'<\frac{2+\alpha-2\gamma_\beta}{2+\alpha+2\gamma_\beta},
\]
the interpolation condition holds for all sufficiently large $p$. Since
\[
 \frac{2+\alpha-2\gamma_\beta}{2+\alpha+2\gamma_\beta}-\beta
 =\frac{\alpha(1-\beta)}{2+\alpha+2\gamma_\beta}
 \ge\frac{\alpha(1-\beta)}{4+\alpha},
\]
we may take
\[
 1-\beta_{j+1}=\frac{8+\alpha}{8+2\alpha}(1-\beta_j).
\]
Starting with $\beta_0=0$ from Proposition~\ref{prop:uniform-gradient},
we reach every fixed target below one in finitely many steps.
Since only finitely many strict interpolation inequalities are
required, one sufficiently large finite $p$ satisfies all of them.

The estimates are uniform for each prescribed class of positive
$C^\alpha$ densities. Only finitely many enlargements of the density
bounds and reductions of the time interval are needed. By comparison and uniform convergence,
passing $\varepsilon\downarrow0$ proves \eqref{eq:para-all-gradient}.
\end{proof}

\Needspace{5\baselineskip}
The gradient bounds also give H\"older estimates for
$\partial_\nu u$ on the boundary and $(u-\Phi)/\rho$ on $\bar\Omega$.

\begin{corollary}
\label{prop:boundary-normal-slope-holder}
Under the hypotheses of Proposition~\ref{prop:para-all-gradient},
for every fixed $0<\beta<1$,
\begin{equation}
 \|\gamma\|_{C^\beta(\bar\Omega)}
 +\|\partial_\nu u\|_{C^\beta(\partial\Omega)}\le C_\beta.
 \label{eq:global-quotient-holder}
\end{equation}
Moreover,
\begin{equation}
 |u(x)-\varphi(p)-Du(p)\cdot(x-p)|
 \le C_\beta|x-p|^{1+\beta},
 \qquad x\in\bar\Omega,\quad p\in\partial\Omega.
 \label{eq:boundary-first-jet}
\end{equation}
\end{corollary}
\begin{proof}
In fixed normal coordinates, we have
\[
 \gamma(y,t)=\frac{t}{\rho(y,t)}
       \int_0^1\partial_t(u-\Phi)(y,st)\,ds.
\]
By Proposition~\ref{prop:para-all-gradient}, $\gamma$ has a uniform
$C^\beta$ bound in a collar, since $t/\rho$ has a uniformly bounded
$C^1$ norm. Away from the boundary, $|\rho|\ge c>0$.
The normal derivative estimate follows from $\partial_\nu u=Du\cdot\nu$. Taylor's formula for a uniformly bounded
$C^{1,\beta}$ extension of $u$ gives \eqref{eq:boundary-first-jet}.
\end{proof}

\subsection{Boundary second derivative estimates}\label{sec:boundary-hessian}

We will estimate the boundary Hessian using the global
$C^{1,\beta}$ bounds, Green function estimates and comparison
arguments. The tangential derivatives of $f$ are controlled after
Neumann heat regularization by Young's inequality for Stieltjes
integrals \cite[Section~10, inequality~(10.9)]{Young}, using only
H\"older norms. Integration by parts in tangential directions produces no boundary term.

\begin{lemma}\label{lem:tangential-young}
Let $m\ge2$ be an integer and let $\Omega\Subset\mathbb R^{m}$ have
smooth boundary. Denote its Neumann heat semigroup by $P_t$ and
its kernel by $p_t$.
Suppose $0<\alpha,\beta<1$ and $\alpha+\beta>1$.
In a fixed open chart neighborhood $Q$, write
$\Omega\cap Q=\{x<\psi(\vartheta)\}\cap Q$ and let
$X=\partial_{\vartheta_j}+\psi_{\vartheta_j}\partial_x$,
where $1\le j\le m-1$.
For real-valued functions $b,f\in C^\infty(\overline{\Omega\cap Q})$ and
$\eta\in C_c^\infty(Q)$, we have
\begin{equation}
 \|P_t(\eta bXf)\|_{L^\infty}\le
 C t^{(\alpha-1)/2}
 \|f\|_{C^\alpha(\overline{\Omega\cap Q})}
 \|b\|_{C^\beta(\overline{\Omega\cap Q})},
 \qquad 0<t\le1. \label{eq:BB1}
\end{equation}
The constant is independent of higher derivatives of $b$ and $f$.

\end{lemma}

\begin{proof}
Set $r=x-\psi(\vartheta)$, where
$\vartheta\in\mathbb R^{m-1}$.
We extend $\psi$ smoothly from a larger coordinate neighborhood
with a uniform $C^3$ bound. Then the map
$\Phi_{\rm g}(\vartheta,r)=(\vartheta,r+\psi(\vartheta))$
is globally bi-Lipschitz, with $\det D\Phi_{\rm g}=1$.
On the support of the chart cutoff, the boundary is given by $r=0$.
The chain rule gives $(Xf)\circ\Phi_{\rm g}
=\partial_{\vartheta_j}(f\circ\Phi_{\rm g})$.
Let $\ell=\sqrt t$ and $y=(\vartheta,r)$. We continue to write
$b,f$ for $b\circ\Phi_{\rm g},f\circ\Phi_{\rm g}$, respectively.
For $\zeta\in\Omega$ and $y\in\Phi_{\rm g}^{-1}(\Omega\cap Q)$, define
\[
 K_t(\zeta;y)=p_t(\zeta,\Phi_{\rm g}(y))\,
                          \eta(\Phi_{\rm g}(y)).
\]
Using \eqref{eq:uniform-neumann-kernel} in real dimension $m$ and
the bi-Lipschitz bounds for $\Phi_{\rm g}$, we obtain, for $s=\vartheta_j$,
\begin{equation}
 \begin{gathered}
 |K_t(\zeta;y)|\le C\ell^{-m}
       e^{-|y-\Phi_{\rm g}^{-1}(\zeta)|^2/(C\ell^2)};\\
 |\partial_sK_t(\zeta;y)|\le C\ell^{-m-1}
       e^{-|y-\Phi_{\rm g}^{-1}(\zeta)|^2/(C\ell^2)}.
 \end{gathered}
 \label{eq:BBkernel}
\end{equation}%
The kernel $K_t(\zeta;\cdot)$ is extended by zero to
$\mathbb R^{m-1}\times(-\infty,0)$. Since the cutoff vanishes near the artificial boundary of the chart,
the zero extension is $C^1$ in $s$. The derivatives of the cutoff
contribute at most $C\ell^{-m}\le C\ell^{-m-1}$ for $\ell\le1$. We extend $b$ and $f$ to the half-space by the infimum
extension formula used in Lemma~\ref{lem:smooth-data}, with exponents
$\beta$ and $\alpha$, respectively. Truncation to their original ranges
preserves their H\"older and supremum bounds. On each $s$-line the
integral against $df$ agrees with the original integral, since the
kernel is supported where $f$ agrees locally with the smooth data.

The $s$-axis admits a partition of unity with bounded overlap, interval
lengths comparable to $\ell$, and derivatives bounded by $C\ell^{-1}$.
Let $\chi_I$ be the function in this partition of unity supported in
an interval $I$, and set $K_I=\chi_I K_t$. Choose $s_I\in I$.
The term $b(s_I)\int_I K_I\,df$ is estimated by integration by parts.
Since $\alpha+\beta>1$, Young's inequality for Stieltjes integrals
\cite[Section~10, inequality~(10.9)]{Young} gives
\[
 \left|\int_I a\,df-a(s_I)(f(\sup I)-f(\inf I))\right|
 \le C_{\alpha,\beta}[a]_{C^\beta(I)}[f]_{C^\alpha(I)}
 |I|^{\alpha+\beta}.
\]
For $a=K_I(b-b(s_I))$, we have
$a(s_I)=a(\inf I)=a(\sup I)=0$, since $\chi_I$ is supported
strictly inside $I$. The correction term therefore vanishes.
The product inequality for H\"older seminorms and the Lipschitz
bound for $K_I$ give
$[a]_{C^\beta}\le C[b]_{C^\beta}
(\|K_I\|_{L^\infty(I)}+\ell\|K_I'\|_{L^\infty(I)})$.
For each fixed choice of the transverse coordinates, we obtain
\[
 \begin{split}
 \left|\int_I K_I b\,df\right|
 \le C[f]_{C^\alpha}\bigl\{
 &\|b\|_{L^\infty(I)}\ell^{\alpha+1}\|K_I'\|_{L^\infty(I)}\\
 &+[b]_{C^\beta}\ell^{\alpha+\beta}
 (\|K_I\|_{L^\infty(I)}+\ell\|K_I'\|_{L^\infty(I)})\bigr\}.
 \end{split}
\]
Using \eqref{eq:BBkernel} and the bounded overlap of the intervals,
we sum the Gaussian factors. Integration in the remaining $m-1$
variables gives a factor at most $C\ell^{m-1}$, even with $r<0$.
It follows that
\[
 |P_t(\eta bXf)(\zeta)|\le C[f]_{C^\alpha}\{
 \|b\|_{L^\infty(\Omega\cap Q)}\ell^{\alpha-1}
 +[b]_{C^\beta}\ell^{\alpha+\beta-1}\}.
\]%
The restriction $r<0$ involves only a transverse variable. Hence no
boundary term arises in the $s$-integrals, proving
\eqref{eq:BB1}.

\end{proof}

We next use the global $C^{1,\beta}$ bound to estimate the Green
function as its pole approaches the boundary.

\begin{lemma}\label{lem:green-pole-decay}
Let $0<\beta<1$ and let $u\in C^\infty(\bar\Omega)$ be strictly
plurisubharmonic. Suppose that
\[
 \lambda\le f=\det(u_{i\bar j})\le\Lambda,\qquad
 \|u\|_{C^{1,\beta}(\bar\Omega)}\le B_\beta.
\]
For smooth $h$, define the Dirichlet solution operator $R_u$ by
\[
 fL_uR_uh=h\quad\text{in }\Omega,\qquad
 R_uh=0\quad\text{on }\partial\Omega.
\]
For $h\in L^\infty(\Omega)$, set
$R_uh(x)=-\int_\Omega G_x(y)h(y)\,dy$, with the normalization
\eqref{eq:FG1}.
For every $0<\gamma<\beta$,
\begin{equation}
 G_x\in W^{1,1}_0(\Omega),\qquad
 \|DG_x\|_{L^1(\Omega)}\le C_\gamma d(x)^\gamma.
 \label{eq:green-pole-decay}
\end{equation}
Consequently, for bounded $k$ and $0<t\le1$,
\begin{equation}
 |R_u(I-P_t)k(x)|\le C_\gamma\sqrt t\,d(x)^\gamma
                                      \|k\|_{L^\infty}.
 \label{eq:green-pole-heat}
\end{equation}
The constants depend only on the fixed geometry,
$n,\lambda,\Lambda,\beta,B_\beta$ and the chosen $\gamma$.
\end{lemma}

\begin{proof}
For the defining function in \eqref{eq:CG2}, we have $-\rho\asymp d$.
After subtracting a constant vector from $Du$, we apply the divergence
theorem to interior balls and boundary caps. Using the $C^{1,\beta}$
bound, we obtain
\[
 \int_{B_r(z)\cap\Omega}\Delta u\,dy\le Cr^{2n-1+\beta}.
\]
Fix $0<a<\beta$ and set $d\nu_a=(-\rho)^{-a}\Delta u\,dy$.
This finite measure, extended by zero, satisfies
\[
 \nu_a(B_r(z))\le C_a r^{s_a},\qquad
                   s_a=2n-1+\beta-a.
\]
Away from the boundary, this follows from $-\rho\ge cr$.
In a boundary ball, each layer $-\rho\asymp2^{-j}r$ is covered by
at most $C2^{j(2n-1)}$ balls of radius $2^{-j}r$. By the preceding estimate for $\Delta u$,
its weighted mass is at most
$Cr^{2n-1}(2^{-j}r)^{\beta-a}$. The sum converges because $a<\beta$.

Let $E\subset\Omega$ have finite perimeter in $\mathbb R^{2n}$.
At each Lebesgue density-one point $x$ of $E$, choose $r_x>0$ such that
$|E\cap B_{r_x}(x)|=|B_{r_x}(x)|/2$.
Such radii exist by continuity of the volume ratio, which tends to $1$
as $r\downarrow0$ and to $0$ as $r\to\infty$. They are uniformly bounded.
Select pairwise disjoint balls $B_i=B_{r_i}(x_i)$ from this family, taking
each radius greater than half the supremum of the remaining radii and
discarding intersecting balls \cite{Gustin}.
In an infinite selection, boundedness and disjointness force $r_i\to0$.
Thus every ball in the original family is eventually selected or discarded.
Every discarded ball lies in the corresponding $B_{5r_i}(x_i)$,
so these enlarged balls cover all density-one points of $E$.
The relative isoperimetric inequality \cite[Section~5.6]{EG} gives
\[
 \sum_i r_i^{2n-1}
 \le C\sum_i\operatorname{Per}(E;B_i)
 \le C\operatorname{Per}(E),
\]
where $\operatorname{Per}(E;B_i)$ is the perimeter of $E$ inside $B_i$.
Since $\nu_a$ is absolutely continuous and $s_a>2n-1$, the growth bound gives
\[
 \begin{aligned}
 \nu_a(E)^{(2n-1)/s_a}
 &\le C_a\left(\sum_i r_i^{s_a}\right)^{(2n-1)/s_a}\\
 &\le C_a\sum_i r_i^{2n-1}
 \le C_a\operatorname{Per}(E).
 \end{aligned}
\]%
For the zero extension of $b\in W^{1,1}_0(\Omega)$, we use
$|b(x)|=\int_0^\infty\mathbf1_{\{|b|>t\}}(x)\,dt$.
By Minkowski's integral inequality \cite[Section~7.1]{GT} and the
coarea formula \cite[Chapter~5]{EG}, we obtain
\begin{equation}
 \|b\|_{L^{s_a/(2n-1)}(\nu_a)}
 \le\int_0^\infty\nu_a(\{|b|>t\})^{(2n-1)/s_a}\,dt
 \le C_a\|Db\|_{L^1}.
 \label{eq:weighted-w11-trace}
\end{equation}%
Localization on cubes yields
\begin{equation*}
 \|b\|_{L^\theta(\nu_a)}\le C\|Db\|_{L^1}^{1-e}\|b\|_{L^1}^e,
 \quad e=\frac{s_a}{\theta}-(2n-1),
 \quad 1<\theta<\frac{s_a}{2n-1}.
\end{equation*}
Here $0<e<\beta-a$. Let $\{Q\}$ be a partition of
$\mathbb R^{2n}$ into cubes of side $\ell$. By
\eqref{eq:weighted-w11-trace} with standard cutoffs for $Q\subset2Q$,
H\"older's inequality and $\nu_a(Q)\le C\ell^{s_a}$, we have
\[
 \left(\int_Q|b|^\theta\,d\nu_a\right)^{1/\theta}
 \le C\ell^e\bigl(\|Db\|_{L^1(2Q)}
                 +\ell^{-1}\|b\|_{L^1(2Q)}\bigr).
\]
Summing over $Q$ and using bounded overlap and
$\|\cdot\|_{\ell^\theta}\le\|\cdot\|_{\ell^1}$, we obtain
\[
 \|b\|_{L^\theta(\nu_a)}
 \le C\ell^e\bigl(\|Db\|_{L^1}+\ell^{-1}\|b\|_{L^1}\bigr).
\]
Taking $\ell=\|b\|_{L^1}/\|Db\|_{L^1}$ for $b\ne0$ gives the desired estimate.

Let $A$ be the real symmetric matrix defined by \eqref{eq:real-cofactor}.
For $N>0$, put $G_N=\min(G_x,N)$ and
$\mu_N=-\operatorname{div}(A DG_N)$.
The measure $\mu_N$ is positive and supported on $\{G_x=N\}$.
For $\varepsilon>0$, we have
$\operatorname{div}(A D[(\varepsilon-\rho)^a])\le0$. Since
$(G_x-N)_+$ has compact support in $\Omega$, the identity
\[
 \mu_N=\delta_x+\operatorname{div}(A D(G_x-N)_+)
\]
gives
\[
 \begin{aligned}
 \int_\Omega(\varepsilon-\rho)^a\,d\mu_N
 &= (\varepsilon-\rho(x))^a
    +\int_\Omega (G_x-N)_+\operatorname{div}(A D[(\varepsilon-\rho)^a])\,dy\\
 &\le(\varepsilon-\rho(x))^a.
 \end{aligned}
\]
For $c>0$, let $\Psi(0)=0$, $\Psi'(s)=(c+s)^{-\theta}$,
and $F(0)=0$, $F'=\Psi$.
Testing with $(\varepsilon-\rho)^a\Psi(G_N)$ gives the cross term
\[
 \int_\Omega A D F(G_N)\cdot D[(\varepsilon-\rho)^a]\,dy
 =-\int_\Omega F(G_N)\operatorname{div}(A D[(\varepsilon-\rho)^a])\,dy
 \ge0.
\]
There is no boundary term because $F(G_N)=0$ there. Since
$\Psi\le c^{1-\theta}/(\theta-1)$, letting
$\varepsilon\downarrow0$ gives
\[
 \int_\Omega\frac{(-\rho)^a DG_N\cdot A DG_N}
                         {(c+G_N)^\theta}\,dy
 \le \frac{c^{1-\theta}}{\theta-1}(-\rho(x))^a.
\]
For each fixed smooth operator, these tests are justified with the
regularized weight before passing to the limit, so no singular
boundary term occurs.

By the Cauchy--Schwarz inequality, \eqref{eq:FG2} and the
interpolation estimate above,
\[
 \|DG_N\|_{L^1}^2\le C(-\rho(x))^a c^{1-\theta}
 \bigl(c+\|DG_N\|_{L^1}^{1-e}\|G_N\|_{L^1}^{e}\bigr)^\theta .
\]
Take $c=\|DG_N\|_{L^1}^{1-e}\|G_N\|_{L^1}^{e}$. Comparison with the
defining function gives $\|G_N\|_{L^1}\le-R_u1(x)\le C(-\rho(x))$.
Hence
\[
 \|DG_N\|_{L^1}
 \le C(-\rho(x))^{a/(1+e)}\|G_N\|_{L^1}^{e/(1+e)}
 \le C(-\rho(x))^{(a+e)/(1+e)}.
\]
Choosing $\gamma<a<\beta$ and $0<e<\beta-a$ gives
$\|DG_N\|_{L^1}\le C_\gamma d(x)^\gamma$. Since
$DG_N=\mathbf1_{\{G_x<N\}}DG_x$ almost everywhere, the bound and
$\|G_N-G_x\|_{L^1}\to0$ give convergence in $W^{1,1}_0(\Omega)$.
Applying the estimate following \eqref{eq:green-bilinear} to $b=G_x$
gives \eqref{eq:green-pole-heat}.
\end{proof}

\begin{lemma}\label{lem:resolvent-decomposition}
Assume the hypotheses of Lemma~\ref{lem:green-pole-decay}.
For $k\in L^\infty(\Omega)$,
\begin{equation}
 \begin{gathered}
 |R_uk(x)|\le Cd(x)\|k\|_{L^\infty},\\
 \|R_uk\|_{L^\infty}
 +\|\partial_\nu R_uk\|_{L^\infty(\partial\Omega)}
 \le C\|k\|_{L^\infty}.
 \end{gathered}
 \label{eq:BB2b}
\end{equation}
For $h\in C^\infty(\bar\Omega)$ and $t_j=4^{-j}$, set $r_j=2^{-j}$.
Then
\begin{equation}
 h=P_1h+\sum_{j\ge1}(I-P_{3t_j})P_{t_j}h.
 \label{eq:BB3}
\end{equation}
For every $0<\gamma<\beta$,
\begin{equation}
 |R_uh(x)|\le Cd(x)\|P_1h\|_{L^\infty}
 +C_\gamma d(x)^\gamma
       \sum_{j\ge1}r_j\|P_{t_j}h\|_{L^\infty}.
 \label{eq:dyadic-potential-source}
\end{equation}
If, for some $0<\sigma<1$,
\[
 \|\partial_\nu R_u(I-P_t)k\|_{L^\infty(\partial\Omega)}
 \le C_\sigma t^{\sigma/2}\|k\|_{L^\infty},\qquad 0<t\le1,
\]
then
\begin{equation}
 \|\partial_\nu R_uh\|_{L^\infty(\partial\Omega)}
 \le C\|P_1h\|_{L^\infty}
 +C_\sigma\sum_{j\ge1}r_j^\sigma\|P_{t_j}h\|_{L^\infty}.
 \label{eq:dyadic-normal-source}
\end{equation}
The constants depend on the data in Lemma~\ref{lem:green-pole-decay}.
The constant in \eqref{eq:dyadic-normal-source} also depends on the
assumed estimate for the normal derivative.
\end{lemma}

\begin{proof}
The global strictly PSH defining function satisfies $fL_u\rho\ge c>0$.
By comparison with $\pm C\|k\|_{L^\infty}(-\rho)$, we obtain the
pointwise bound. For each fixed operator, the $W^{2,p}$
Dirichlet estimate with $p>2n$ and Sobolev embedding
\cite[Theorems~9.15 and~7.26]{GT} give $R_uk\in C^1(\bar\Omega)$.
Boundary differentiation then gives \eqref{eq:BB2b}.
By the semigroup property,
\[
 P_1h+\sum_{j=1}^J(I-P_{3t_j})P_{t_j}h=P_{t_J}h.
\]%
Using \eqref{eq:green-pole-heat} on each summand gives
\eqref{eq:dyadic-potential-source}. The assumed estimate for the boundary normal derivative
gives \eqref{eq:dyadic-normal-source}, and the term involving $P_1h$
is controlled by \eqref{eq:BB2b}. For each fixed smooth $h$,
\eqref{eq:neumann-heat-bounds} gives $P_{t_J}h\to h$ uniformly.
No Neumann boundary condition is needed for $h$. By \eqref{eq:BB2b},
$R_uP_{t_J}h\to R_uh$ uniformly, and their boundary normal derivatives
also converge uniformly.
\end{proof}

The next lemma bounds the boundary normal derivative by comparison
with $R_uh$. The constants do not depend on derivatives
of the right-hand side.

\begin{lemma}\label{lem:boundary-source-comparison}
Assume the hypotheses of Lemma~\ref{lem:green-pole-decay}.
Let $p\in\partial\Omega$, $D=\Omega\cap B_r(p)$ lie in a fixed boundary
chart. Suppose that
$\Psi\in C^\infty(\bar D)$ and $h\in C^\infty(\bar\Omega)$ are real-valued, and
\[
 \Psi=0\quad\text{on }\partial\Omega\cap B_r(p),\qquad
 L_u\Psi\le h/f+C_0\sum_i u^{i\bar i}\quad\text{in }D,\qquad C_0\ge0.
\]
With $Z=R_uh$,
\begin{equation}
 \partial_\nu\Psi(p)\le\partial_\nu Z(p)
 +C\left[C_0+r^{-2}
       \bigl(\|Z\|_{L^\infty(D)}+\|\Psi\|_{L^\infty(D)}\bigr)\right].
 \label{eq:boundary-source-normal}
\end{equation}
Here $C$ depends only on the fixed defining function and geometry.
\end{lemma}

\begin{proof}
By \eqref{eq:CG2}, $L_u\rho\ge c_0\sum_i u^{i\bar i}$. Set
\[
 M=r^{-2}\bigl(\|Z\|_{L^\infty(D)}+\|\Psi\|_{L^\infty(D)}\bigr),
 \quad A=(M+C_0)/c_0,\quad \widetilde{Z}(x)=Z(x)+A\rho(x)-M|x-p|^2.
\]
On $\partial D$, we have $\widetilde{Z}\le\Psi$, while
$L_u\widetilde{Z}\ge h/f+C_0\sum_i u^{i\bar i}\ge L_u\Psi$ in $D$.
The maximum principle gives $\widetilde{Z}\le\Psi$ throughout $D$.
Since $\widetilde{Z}(p)=\Psi(p)=0$, outward normal differentiation gives
\eqref{eq:boundary-source-normal}.
\end{proof}

Let $\Omega\Subset\mathbb C^n$ be a fixed smooth bounded strongly
pseudoconvex domain. Fix a smooth real-valued boundary value $\varphi$ and
$0<\lambda\le\Lambda$. Let $u\in C^\infty(\bar\Omega)$ be a
strictly PSH solution satisfying
\begin{equation}
 \det(u_{i\bar j})=f,\qquad \lambda\le f\le\Lambda,
 \qquad u|_{\partial\Omega}=\varphi,
 \qquad\|u\|_{C^{1,\beta}(\bar\Omega)}\le B,
 \quad0<\beta<1.                                     \label{eq:MC1}
\end{equation}
All constants below depend only on
$n,\Omega,\varphi,\lambda,\Lambda,\beta,B$.
Here the dependence on $\Omega,\varphi$ is through the
quantitative $C^3$ geometric and boundary data fixed in
Section~\ref{sec:preliminaries}, including the continuity moduli of the
third derivatives of the boundary charts and boundary extension.
In particular, they are independent of interior second derivatives
of $u$ and derivatives of $f$.

The next lemma combines the barrier inequality in
\cite[Lemma~2.1]{Guan} with a gradient bound of order $R^\beta$
for the comparison function.

\begin{lemma}\label{lem:matching-barrier}
Let $\rho,\Phi$ be the fixed functions in
Lemma~\ref{lem:geometric-barriers}. Under \eqref{eq:MC1}, there are
positive constants $K_0,C_1,c,C$ and a radius
$R_0>0$ with the following property.  For every $p\in\partial\Omega$ and $0<R\le R_0$, set
\[
 v_{p,R}=u-\Phi-(\gamma(p)+C_1R^\beta)\rho-K_0\rho^2.
\]
Then, on $\Omega\cap B_R(p)$,
\[
 v_{p,R}\ge0,\qquad |Dv_{p,R}|\le CR^\beta,\qquad
 L_uv_{p,R}\le-c\bigl(1+\sum_i u^{i\bar i}\bigr).
\]
Moreover, $v_{p,R}=0$ on $\partial\Omega\cap B_R(p)$.
The constants and $R_0$ depend only on the data listed after \eqref{eq:MC1}.

\end{lemma}

\begin{proof}
Since $D(u-\Phi)(p)=\gamma(p)D\rho(p)$, the common $C^{1,\beta}$ bound
gives, for $z\in\bar\Omega\cap B_{2R}(p)$,
\[
 |D(u-\Phi)(z)-\gamma(p)D\rho(z)|\le CR^\beta.
\]
Integrating along the inward normal segment from the nearest
boundary point gives
\[
 |(u-\Phi)/\rho-\gamma(p)|\le CR^\beta
 \quad\text{on }\Omega\cap B_R(p).
\]
Also $-\rho\le CR$ there. For sufficiently small fixed $R_0$,
we therefore have $-\rho<\kappa R^\beta$ whenever $R\le R_0$.
Lemma~\ref{lem:cap-comparison} applies with
$E=R^\beta$, $H=\kappa E$ and $D=\Omega\cap B_R(p)$.
Since $E/H=\kappa^{-1}$ is fixed, the comparison function in
Lemma~\ref{lem:cap-comparison} has the stated form after setting
$C_1=D_0$ and increasing the fixed quadratic coefficient $K_0$
by $\kappa^{-1}$.
Consequently $v_{p,R}\ge0$, $v_{p,R}=0$ on
$\partial\Omega\cap B_R(p)$, and
$L_uv_{p,R}\le-c(1+\sum_i u^{i\bar i})$.
Finally,
\[
 |Dv_{p,R}|\le CR^\beta+C_1R^\beta|D\rho|
                     +2K_0|\rho||D\rho|\le CR^\beta,
\]
since $R\le1$. All choices depend only on the prescribed data,
$\beta$ and the common $C^{1,\beta}$ bound.
\end{proof}

The following estimate for the boundary normal derivative requires
no additional $L^p$ bound for the Hessian with $p>1$.

\begin{lemma}\label{lem:matching-normal-heat}
Assume \eqref{eq:MC1}. For every $0<\sigma<\beta/(2-\beta)$ and
$h\in L^\infty(\Omega)$,
\begin{equation}
 \|\partial_\nu R_u(I-P_t)h\|_{L^\infty(\partial\Omega)}
       \le C_\sigma t^{\sigma/2}\|h\|_{L^\infty},
 \qquad 0<t\le1.
 \label{eq:matching-normal-heat}
\end{equation}
The constant $C_\sigma$ depends only on $\sigma$ and the data
listed after \eqref{eq:MC1}.
\end{lemma}

\begin{proof}
By homogeneity, assume $\|h\|_{L^\infty}=1$. Choose $0<\gamma<\beta$
such that $\sigma<\beta/(2-\gamma)$, and put $w=R_u(I-P_t)h$.
For $p\in\partial\Omega$ and $0<R\le R_0$, let $v_{p,R}$ be
given by Lemma~\ref{lem:matching-barrier} and set
$D_R=\Omega\cap B_R(p)$.
By \eqref{eq:green-pole-heat}, $|w|\le C_\gamma\sqrt t\,R^\gamma$
on $D_R$. Set
\[
 B=C_\gamma\sqrt t\,R^{\gamma-2},\qquad
 Y=C_0(1+B)v_{p,R}+B|z-p|^2.
\]
Since $L_u|z-p|^2=\sum_i u^{i\bar i}$, taking
$C_0c\ge\max\{1,2/\lambda\}$ gives
\[
 -L_uY\ge C_0c(1+B)+\bigl(C_0c(1+B)-B\bigr)\sum_i u^{i\bar i}
            \ge 2/\lambda\ge |L_uw|.
\]
On $\Omega\cap\partial B_R(p)$, $B|z-p|^2$ dominates $|w|$.
On $\partial\Omega\cap B_R(p)$, we have $w=0$ and $Y\ge0$.
Comparison therefore gives
$|w|\le Y$ on $D_R$. At $p$, the quadratic term has zero normal
derivative, while $v_{p,R}\ge0$ and $v_{p,R}(p)=0$ imply
$\partial_\nu v_{p,R}(p)\le0$. Outward normal differentiation yields
\[
 |\partial_\nu w(p)|\le-\partial_\nu Y(p)
 \le C_\gamma\{R^\beta+\sqrt t\,R^{\beta+\gamma-2}\}.
\]
Taking $R=t^{1/[2(2-\gamma)]}$ for small $t$, we obtain
\[
 |\partial_\nu w(p)|\le C_\gamma t^{\beta/[2(2-\gamma)]}
                         \le C_\sigma t^{\sigma/2}.
\]
For the remaining times, \eqref{eq:BB2b} gives the same estimate.
\end{proof}

\begin{proposition}\label{prop:gradient-bootstrap}
Let $n\ge2$ and $0<\alpha<1$. Fix a bounded smooth strongly pseudoconvex domain
$\Omega$, a smooth boundary value $\varphi$, and constants
$0<\lambda\le\Lambda$ and $F_0<\infty$.
Every smooth strictly plurisubharmonic solution of
\eqref{intro:dirichlet} with $\lambda\le f\le\Lambda$ and
$\|f\|_{C^\alpha(\bar\Omega)}\le F_0$ satisfies
\[
 \sup_{\partial\Omega}|D^2u|\le C,
 \qquad \|u\|_{C^{1,\beta}(\bar\Omega)}\le C_\beta,
\]
with $\beta=\beta_n$ as defined in \eqref{eq:PIA12}.
The constants are uniform under the fixed $C^3$ approximation in
Lemma~\ref{lem:smooth-data}, including the third-derivative
continuity moduli.
\end{proposition}

\begin{proof}
For every prescribed finite $C^\alpha$ bound on positive densities,
Proposition~\ref{prop:para-all-gradient} gives a common global
$C^{1,\beta}$ bound for every fixed $0<\beta<1$. Take
$\beta=\beta_n$ from \eqref{eq:PIA12}. Since $1-\beta_n<\alpha/4$,
we have $\alpha+\beta_n>1$ and
\[
 \frac{\beta_n}{2-\beta_n}-\left(1-\frac\alpha2\right)
 =\frac\alpha2-\frac{2(1-\beta_n)}{2-\beta_n}>0.
\]
We may therefore fix $\sigma=1-\alpha/2$, for which
$\alpha+\sigma-1=\alpha/2>0$.

In a fixed small boundary chart of radius $2r_0$, let $T_j$ be the
projected fields in Lemma~\ref{lem:para-boundary-trace}.
The commutator estimate \eqref{eq:para-nearhol-source} also applies
to $u$.
With $w_j=T_j(u-\Phi)$ and
$\mathcal E=\sum_j|\partial w_j|_u^2$, it gives
\[
 |L_uw_k-T_k\log f|\le C\sum_i u^{i\bar i}+Cr_0\Bigl(\mathcal E\sum_i u^{i\bar i}\Bigr)^{1/2}.
\]
The terms from the fixed holomorphic coordinate change are bounded
by $C\sum_i u^{i\bar i}$.
Using a quadratic correction as in \cite[Section~3.3(c)]{Guan2009},
we set $\Psi_k^\pm=\pm w_k-\sum_jw_j^2$.
Using the common gradient bound and Young's inequality, we absorb
the commutator errors into the term $-2\mathcal E$ and obtain
\[
 L_u\Psi_k^\pm\le
 \frac{\pm T_kf-2\sum_jw_jT_jf}{f}+C\sum_i u^{i\bar i}.
\]
Let $\eta\in C_c^\infty(B_{2r_0})$ equal one on $B_{r_0}$, and set
\[
 h_k^\pm=\eta\left(\pm T_kf-2\sum_jw_jT_jf\right).
\]
In fixed coordinates adapted to the level sets of $\rho$, the fields
$T_j$ have only tangential components. We compose the Neumann heat kernel with the fixed chart and
include the Jacobian determinant. The resulting kernel satisfies the
bounds used in Lemma~\ref{lem:tangential-young}.
The coordinate coefficients of $T_j$ have common $C^1$ bounds.
Since $\|w_j\|_{C^\beta}\le C$
and $\alpha+\beta>1$, that lemma gives
$\|P_th_k^\pm\|_{L^\infty}\le Ct^{(\alpha-1)/2}$ for $0<t\le1$.
By Lemmas~\ref{lem:matching-normal-heat} and
\ref{lem:resolvent-decomposition}, we obtain
\[
 \begin{aligned}
 \|R_uh_k^\pm\|_{L^\infty}
 &\le C\left(1+\sum_{j\ge1}2^{-j\alpha}\right)\le C,\\
 \|\partial_\nu R_uh_k^\pm\|_{L^\infty(\partial\Omega)}
 &\le C\left(1+\sum_{j\ge1}2^{-j(\alpha+\sigma-1)}\right)\le C.
 \end{aligned}
\]
On $\partial\Omega$, $w_j=0$ and
$\partial_\nu\Psi_k^\pm=\pm\partial_\nu w_k$. The gradient bound
gives $\|\Psi_k^\pm\|_{L^\infty}\le C$. By
Lemma~\ref{lem:boundary-source-comparison}, $|\partial_\nu w_k|\le C$.
The real mixed derivatives are bounded because
\[
 \partial_\nu w_k=D^2(u-\Phi)(\nu,T_k)
                  +D(u-\Phi)\cdot D_\nu T_k
\]
and the tangential fields have bounded first derivatives.
Using the prescribed boundary value, \eqref{eq:GC20} and
\eqref{eq:boundary-schur}, we obtain a bound for the full real Hessian on $\partial\Omega$.
\end{proof}

\section{Global Hessian estimates}\label{sec:global-hessian}

In this section, we first derive integral estimates for the Hessian. Together
with the boundary Hessian bound, they yield a global bound by
comparison with solutions whose densities are locally constant.

\subsection{Interior Hessian and integral estimates}

The next lemma makes the dependence on the $C^{1,\beta}$ norm
explicit in the comparison argument of Li--Li--Zhang
\cite[Theorem~1.2]{LLZ}. We also use the derivative estimate of
\cite{RS}, in the form given in
\cite[Proposition~2.4, equation~(2.9)]{LLZ}.

\begin{lemma}\label{lem:polynomial-interior}
Let $n\ge2$, $0<\alpha<1$, $0<\lambda\le\Lambda$ and $F_0<\infty$
be fixed. Suppose
\begin{equation}
 1-\frac\alpha{16(n+1)}\le\beta<1.
 \label{eq:PIA1}
\end{equation}
If $u$ is smooth and strictly PSH on $B_2\subset\mathbb C^n$,
\[
 \det(u_{i\bar j})=f,\quad
 \lambda\le f\le\Lambda,\quad
 \|f\|_{C^\alpha(B_2)}\le F_0,\quad
 \|u\|_{C^{1,\beta}(B_2)}\le K,\quad K\ge1,
\]
then
\begin{equation}
 |D^2u(0)|\le CK^{n+1}.
 \label{eq:PIA2}
\end{equation}
The constant depends only on the fixed data, $\beta$, and a fixed
radial mollifier.
\end{lemma}

\begin{proof}
After subtracting the affine Taylor polynomial of $u$ at zero, set
$u^t=\varrho_{t^2}*u$, where $\varrho$ is a fixed radial
nonnegative unit-mass mollifier. For sufficiently small $t>0$,
\[
 u\le u^t\le u+CKt^{2+2\beta},\qquad
 |D^2u^t|\le CKt^{-2(1-\beta)}\quad\hbox{on }B_t.
\]
By \cite{CKNS}, there is a smooth strictly PSH solution
\[
 \det((\phi_t)_{i\bar j})=f(0)\quad\hbox{in }B_t,
 \qquad \phi_t=u^t\quad\hbox{on }\partial B_t.
\]
Concavity of $\det^{1/n}$ and the comparison principle
\cite[Theorem~A]{BT}, applied with the correction
$Ct^\alpha(|z|^2-t^2)$, give
\begin{equation}
 \|\phi_t-u\|_{L^\infty(B_t)}\le CKt^{2+\alpha}.
 \label{eq:PIA4}
\end{equation}
Here $2+2\beta>2+\alpha$ by \eqref{eq:PIA1}, and $K\ge1$.
After quadratic rescaling to $B_1$, the $C^{1,1}$ norm of the
rescaled boundary value is at most $CKt^{-2(1-\beta)}$.
The interior estimate on a ball, linear in the $C^{1,1}$ norm of
the boundary value for constant density
\cite[Proposition~6.6, Theorem~6.7 and Remark~6.8]{BT}, gives
\begin{equation}
 \|D^2\phi_t\|_{L^\infty(B_{t/2})}
       \le CKt^{-2(1-\beta)}.
 \label{eq:PIA5}
\end{equation}
Applying
\cite[Proposition~2.4, equation~(2.9)]{LLZ} to $\phi_t$ on $B_{t/2}$, we obtain
\begin{equation}
 t\|D((\phi_t)_{i\bar j})\|_{L^\infty(B_{t/4})}
 \le\frac{C(n)}{f(0)}\|\Delta\phi_t\|_{L^\infty(B_{t/2})}^{n+1}.
 \label{eq:calabi-estimate}
\end{equation}
Consequently
\[
 t\|\nabla\Delta\phi_t\|_{L^\infty(B_{t/4})}
 +\|\Delta\phi_t\|_{L^\infty(B_{t/4})}
 \le CK^{n+1}t^{-2(1-\beta)(n+1)}.
\]

Fix one sufficiently small $t_0>0$ and let $t_k=2^{-k}t_0$.
For $w_k=\phi_{t_{k-1}}-\phi_{t_k}$ and $t=t_k$,
\eqref{eq:PIA4} and the Schauder estimate for the Poisson equation
\cite[Theorem~4.8]{GT} give
\[
 \begin{gathered}
 \|w_k(t\,\cdot)\|_{L^\infty(B_{1/8})}\le CKt^{2+\alpha};\\
 \|w_k(t\,\cdot)\|_{C^{2,1/2}(B_{1/8})}
       \le CK^{n+1}t^{2-2(1-\beta)(n+1)}.
 \end{gathered}
\]
Using the interpolation inequality \cite[Lemma~6.32]{GT}
and rescaling, we obtain
\[
 \begin{aligned}
 \|D^2w_k\|_{L^\infty(B_{t/16})}
 &\le Ct^{-2}(Kt^{2+\alpha})^{1/5}
       (K^{n+1}t^{2-2(1-\beta)(n+1)})^{4/5}+CKt^\alpha\\
 &\le CK^{(4n+5)/5}t^{[\alpha-8(n+1)(1-\beta)]/5}
       +CKt^\alpha\\
 &\le CK^{n+1}t^{\alpha/10}.
\end{aligned}
\]
Indeed, $(4n+5)/5\le n+1$ and
$[\alpha-8(n+1)(1-\beta)]/5\ge\alpha/10$ by \eqref{eq:PIA1}.
These bounds are summable in $k$. Using \eqref{eq:PIA5} at the fixed radius
$t_0$, we obtain
$\|D^2\phi_{t_k}\|_{L^\infty(B_{t_k/16})}\le CK^{n+1}$.
For a real unit vector $\xi$ and $h=t_k/32$, \eqref{eq:PIA4} gives
\[
 \frac{|u(h\xi)+u(-h\xi)-2u(0)|}{h^2}
       \le CK^{n+1}+CKt_k^\alpha.
\]
Letting $k\to\infty$ proves \eqref{eq:PIA2}, since $u$ is smooth.
\end{proof}

For $p\ge1$ and $0\le\kappa<1/p$, the bound
$|D^2u|\le C(1+d^{-\kappa})$ gives, by integration in the fixed
boundary collar,
\begin{equation}
 \int_{\Omega\cap B_R(z)}|D^2u|^p\le C_pR^{2n-p\kappa},
 \qquad z\in\partial\Omega,\quad0<R\le R_0.
 \label{eq:local-polynomial-moment}
\end{equation}
Indeed, the boundary projection has area at most $CR^{2n-1}$, and
$\int_0^{CR}t^{-p\kappa}\,dt\le C_pR^{1-p\kappa}$.

\begin{corollary}\label{cor:finite-moment}
Let $0<\alpha<1$. Consider a family of smooth Dirichlet
solutions whose densities have common positive lower and upper bounds
and a common $C^\alpha$ bound. Suppose the solutions have a common
global $C^{1,\beta_n}$ bound, where
\begin{equation}
 \beta_n=1-\frac{\alpha^2}{64n^2(n+1)},\qquad
 p_n=2n^2,\qquad
 \kappa_n=(n+1)(1-\beta_n)=\frac{\alpha^2}{64n^2}.
 \label{eq:PIA12}
\end{equation}
Then
\begin{equation}
 |D^2u(x)|\le C(1+d(x)^{-\kappa_n}),\qquad p_n\kappa_n<1,
 \label{eq:PIA11}
\end{equation}
and
\begin{equation}
 \|u\|_{W^{2,p_n}(\Omega)}\le C,\qquad
 \int_{\Omega\cap B_R(z)}|D^2u|^{p_n}
       \le CR^{2n-p_n\kappa_n}\quad(z\in\partial\Omega,\ 0<R\le R_0).
 \label{eq:PIA13}
\end{equation}
The constants depend only on the common bounds above and the geometry
of the fixed boundary collar.
\end{corollary}

\begin{proof}
The hypotheses of Lemma~\ref{lem:polynomial-interior} hold, since
\[
 16(n+1)(1-\beta_n)=\frac{\alpha^2}{4n^2}<\alpha.
\]
At $x\in\Omega$, we subtract the affine Taylor polynomial at $x$
and rescale quadratically on a ball of radius comparable to $d(x)$.
The rescaled function has $C^{1,\beta_n}$ norm at most
$Cd(x)^{\beta_n-1}$, and the rescaled density has a uniformly
bounded $C^\alpha$ seminorm.
Thus \eqref{eq:PIA2} gives the pointwise bound in \eqref{eq:PIA11}. Moreover,
\[
 p_n\kappa_n=\frac{\alpha^2}{32}<1.
\]
Equation~\eqref{eq:local-polynomial-moment} and integration over the
fixed collar give \eqref{eq:PIA13}.
\end{proof}

\subsection{The global Hessian estimate}

We will deduce the global Hessian bound from uniform
$C^{1,\beta_n}$ and boundary Hessian bounds.

Throughout this subsection, $n\ge2$, $0<\alpha<1$ and
$\beta=\beta_n$ as in \eqref{eq:PIA12}.

Let $\Omega\Subset\mathbb C^n$ be a fixed smooth bounded strongly
pseudoconvex domain, with fixed smooth boundary data $\varphi$.
We fix $0<\lambda\le\Lambda$ and $F_0<\infty$.
For each smooth positive right-hand side, \cite{CKNS} gives a
smooth strictly PSH solution of \eqref{intro:dirichlet}.
Let $\chi\in C_c^\infty(B_1)$ satisfy $0\le\chi\le1$ and
$\chi=1$ on $B_{1/2}$. Set
$F_1=(3+2\|D\chi\|_{L^\infty(B_1)})F_0$.

\begin{proposition}\label{prop:finite-closure}
Suppose that every smooth right-hand side $g$ satisfying
\begin{equation}
 \lambda\le g\le\Lambda,\qquad
 \|g\|_{C^\alpha(\bar\Omega)}\le F_1                  \label{eq:FC2}
\end{equation}
has a smooth strictly plurisubharmonic solution $u_g$ of
\eqref{intro:dirichlet} with right-hand side $g$ such that
\begin{equation}
 \|u_g\|_{C^{1,\beta}(\bar\Omega)}\le B,
 \qquad\sup_{\partial\Omega}|D^2u_g|\le K, \label{eq:FC3}
\end{equation}
where $B,K>0$ are independent of $g$.
Then every smooth right-hand side $f$ with
$\lambda\le f\le\Lambda$ and $\|f\|_{C^\alpha(\bar\Omega)}\le F_0$
has a solution $u_f$ of \eqref{intro:dirichlet} satisfying
\begin{equation}
 \|D^2u_f\|_{L^\infty}\le C.      \label{eq:FC5}
\end{equation}
The constant $C$ depends only on
$n,\alpha,\lambda,\Lambda,F_0,B,K$
and the fixed geometric and boundary data.
\end{proposition}

\begin{proof}
Choose
$r_0\le\min\{1,R_0\}$ so that the boundary neighborhoods below
lie in the fixed boundary collar and the nearest-point projection
onto $\partial\Omega$ is well defined within distance $r_0$. Here $R_0$ is
the common radius in Lemma~\ref{lem:matching-barrier}.
For $p\in\partial\Omega$ and $0<r<r_0$, set
\begin{equation}
 g_{p,r}(z)=f(z)+\chi((z-p)/r)\{f(p)-f(z)\},
 \qquad \widetilde u=u_{g_{p,r}}.                         \label{eq:FC6}
\end{equation}
Since $g_{p,r}$ is a convex combination of $f(z)$ and $f(p)$,
we have $\lambda\le g_{p,r}\le\Lambda$. The cutoff bound and
$|f(z)-f(p)|\le F_0r^\alpha$ on $B_r(p)$ give
\[
 [g_{p,r}]_\alpha
 \le(2+\|D\chi\|_{L^\infty(B_1)})[f]_\alpha,\qquad
 \|g_{p,r}\|_{C^\alpha}\le F_1.
\]
Thus \eqref{eq:FC3} holds for $\widetilde u$ with common constants.
Since $\chi=1$ on $B_{1/2}$, we have
\begin{equation}
 \det(\widetilde u_{i\bar j})=f(p)
 \quad\text{in }\Omega\cap B_{r/2}(p).
 \label{eq:FC-frozen-determinant}
\end{equation}

Set
\[
 Q_{p,r}(z)=r^2\log\frac{|z-p|^2+r^2}
                   {(\operatorname{diam}\Omega+2)^2+r^2}\le0.
\]%
Its complex Hessian has eigenvalues $r^2/(|z-p|^2+r^2)$ with
multiplicity $n-1$ and $r^4/(|z-p|^2+r^2)^2$ in the radial
direction. It follows that
\[
 \det((Q_{p,r})_{i\bar j})
   =\frac{r^{2n+2}}{(|z-p|^2+r^2)^{n+1}}\ge2^{-(n+1)}
       \quad\hbox{on }B_r(p),
 \qquad\|Q_{p,r}\|_{L^\infty}\le Cr^2(1+|\log r|).
\]
The difference $f^{1/n}-g_{p,r}^{1/n}$ is supported in $B_r(p)$
and has absolute value at most $Cr^\alpha$. By Minkowski's determinant inequality \cite[Section~7.8]{HJ},
$\widetilde u+Cr^\alpha Q_{p,r}$ and
$u_f+Cr^\alpha Q_{p,r}$ have determinants at least $f$ and
$g_{p,r}$, respectively. Their boundary values are at most $\varphi$.
Applying the comparison principle \cite[Theorem~A]{BT} in both
directions, we obtain
\begin{equation}
 \|u_f-\widetilde u\|_{L^\infty}
          \le Cr^{2+\alpha}(1+|\log r|).                                               \label{eq:FC7}
\end{equation}

Let $A$ be the real symmetric matrix in \eqref{eq:real-cofactor}
for $\widetilde u$, normalized by
$\operatorname{div}(A\nabla h)=g_{p,r}L_{\widetilde u}h$.
The rows of $A$ are divergence-free by \eqref{eq:cofactor-divergence}.
Let $\lambda_A,\Lambda_A$ be the smallest and largest real
eigenvalues of $A$.
By \eqref{eq:real-cofactor}--\eqref{eq:FG2} and the determinant bounds,
\begin{equation}
 \lambda_A^{-1}\le C\Delta\widetilde u,\qquad
 \Lambda_A\le C(\Delta\widetilde u)^{n-1}.
                                                               \label{eq:FC8}
\end{equation}

For a fixed real unit vector $e$,
differentiating \eqref{eq:FC-frozen-determinant} twice gives
\[
 L_{\widetilde u}\widetilde u_{ee}
  =\operatorname{tr}\bigl(\bigl[(\widetilde u_{i\bar j})^{-1}
                ((\widetilde u_e)_{i\bar j})\bigr]^2\bigr)\ge0.
\]
The right-hand side is nonnegative because the matrix inside
the square is similar to a Hermitian matrix.
Set $w_e=(\widetilde u_{ee}-K-1)_+$.
By smooth convex approximation of the positive part, we have
\begin{equation}
 \operatorname{div}(A\nabla w_e)\ge0
 \quad\text{weakly in }\Omega\cap B_{r/2}(p).
 \label{eq:FC-subsolution}
\end{equation}

The boundary bound in \eqref{eq:FC3} gives
$\widetilde u_{ee}-K-1\le-1$ on $\partial\Omega$.
Thus $w_e$ vanishes near $\partial\Omega$ for each fixed auxiliary
solution. Extending it by zero and taking $A=I$ outside $\Omega$, we retain
\eqref{eq:FC-subsolution} in $B_{r/2}(p)$.
Local boundedness under integrability assumptions on
$\lambda_A^{-1}$ and $\Lambda_A$ is proved in \cite{BS,Trudinger1971}.
Wang--Zhou \cite{WangZhou} established Sobolev inequalities and
regularity estimates for the linearized complex Monge--Amp\`ere equation.
We track the radius in the Moser iteration.

For $R=r/2$, Corollary~\ref{cor:finite-moment} and
\eqref{eq:FC8}, applied on the boundary-centered ball $B_R(p)$, give
\begin{equation}
 \begin{gathered}
 \left(\frac1{|B_R|}\int_{B_R(p)}\lambda_A^{-p_n}\right)^{1/p_n}
 \left(\frac1{|B_R|}\int_{B_R(p)}\Lambda_A^{p_n/(n-1)}
       \right)^{(n-1)/p_n}\le CR^{-n\kappa_n},\\
 \left(\frac1{|B_R|}\int_{B_R(p)}w_e^{p_n}\right)^{1/p_n}
       \le CR^{-\kappa_n}.
 \end{gathered}
 \label{eq:FC-moment-factor}
\end{equation}
Taking $A=I$ outside $\Omega$ adds only $O(R^{2n})$ to the
coefficient integrals.

After rescaling $B_R(p)$ to $B_1$, let $0<\tau<\sigma\le1$ and choose
$\zeta\in C_c^\infty(B_\sigma)$ with $0\le\zeta\le1$,
$\zeta=1$ on $B_\tau$ and $|D\zeta|\le C/(\sigma-\tau)$.
We test with $\zeta^2w_e^{k-1}$.
The $W^{1,2p_n/(p_n+1)}$ Sobolev inequality
\cite[Theorem~7.10]{GT} then yields
\begin{equation}
 \|w_e\|_{L^{np_nk/((n-1)p_n+n)}(B_\tau)}
 \le\left[\frac{Ck^2R^{-n\kappa_n}}{(\sigma-\tau)^2}\right]^{1/k}
      \|w_e\|_{L^{p_nk/(p_n-n+1)}(B_\sigma)},
 \quad k\ge2.
 \label{eq:FC10}
\end{equation}
The exponents
\[
 k_j=(p_n-n+1)
       \left(\frac{n(p_n-n+1)}{(n-1)p_n+n}\right)^j,
 \qquad \sum_{j\ge0}k_j^{-1}=\frac1n,
\]
give, by \eqref{eq:FC-moment-factor},
\begin{equation}
 \sup_{\Omega\cap B_{r/4}(p)}w_e\le Cr^{-2\kappa_n}.
 \label{eq:FC11}
\end{equation}
The constant is independent of $e$.

Let $q\in\partial\Omega\cap B_{r/32}(p)$ and take $v_{q,r/16}$
from Lemma~\ref{lem:matching-barrier}. For sufficiently large fixed
$C,C_1>0$, the function
\[
 \Psi=Cr^{-2-2\kappa_n}\{C_1v_{q,r/16}+|z-q|^2\}
\]
satisfies $L_{\widetilde u}\Psi\le0$ and dominates $w_e$ on
$\partial(\Omega\cap B_{r/16}(q))$ by \eqref{eq:FC11}, since
$B_{r/16}(q)\subset B_{r/4}(p)$. Applying the weak maximum principle
\cite[Theorem~8.1]{GT} and using $|Dv_{q,r/16}|\le Cr^\beta$ and
$d^2\le dr^\beta$, we obtain, for $0\le d<r/16$,
\begin{equation}
 w_e(q-d\nu(q))\le C_*d\,r^{-1-\theta},\qquad
 \theta=1-\beta+2\kappa_n>0.
 \label{eq:FC12}
\end{equation}
The parameters in \eqref{eq:PIA12} give
\begin{equation}
 \theta=\frac{(2n+3)\alpha^2}{64n^2(n+1)}
       <\frac{\alpha^2}{16n^2}<\frac\alpha{32}.
 \label{eq:FC-theta-small}
\end{equation}

Choose $0<c_*\le\min\{1/64,(2C_*)^{-1}\}$. By
\eqref{eq:FC12}, $\widetilde u_{ee}\le K+2$ for every real unit
vector $e$ if the nearest boundary point lies in $B_{r/32}(p)$
and $d\le c_*r^{1+\theta}$.
Since $0\le\widetilde u_{ee}+\widetilde u_{Je,Je}$, the same bound
for $Je$ gives
\[
 -(K+2)\le\widetilde u_{ee}\le K+2.
\]
Thus
\begin{equation}
 \begin{gathered}
 |D^2\widetilde u(q-d\nu(q))|\le C,\\
 q\in\partial\Omega\cap B_{r/32}(p),\qquad 0\le d\le c_*r^{1+\theta}.
 \end{gathered}
 \label{eq:FC13}
\end{equation}

Let $x\in\Omega$ have a nearest boundary point $p$ and a
sufficiently small distance $d=d(x)$ from the boundary. Choose
$c_*r^{1+\theta}=4d$ and rescale at radius $d/8$, so that
\begin{equation}
 r=(4d/c_*)^{1/(1+\theta)},\qquad s=d/8.
 \label{eq:FC-distance-scaling}
\end{equation}
When $d\le(c_*/4)r_0^{1+\theta}$, we have $r\le r_0$.
For $y\in B_{d/2}(x)$ with nearest boundary point $q_y$,
\[
 d(y)\le\frac{3d}{2}<c_*r^{1+\theta},\qquad
 |q_y-p|\le3d\le\frac{3c_*r}{4}<\frac r{32}.
\]
Thus \eqref{eq:FC13} and \eqref{eq:FC-frozen-determinant} give
\[
 |D^2\widetilde u|\le C,\qquad
 \det(\widetilde u_{i\bar j})=f(p)\quad\text{on }B_{d/2}(x).
\]%

On $B_4$, define
\begin{equation}
 \widetilde{U}(y)=\frac1{s^2}
       \{\widetilde u(x+sy)-\widetilde u(x)
                             -sD\widetilde u(x)\cdot y\}.
 \label{eq:FC-reference-scaling}
\end{equation}
Set
\begin{equation}
 U_*(y)=\widetilde U(y)+s^{-2}(u_f-\widetilde u)(x+sy).
 \label{eq:FC-original-scaling}
\end{equation}
By \eqref{eq:FC6} and \eqref{eq:FC-frozen-determinant},
\begin{equation}
 \det(\widetilde{U}_{i\bar j})=f(p),\qquad
 \det((U_*)_{i\bar j})=f(x+sy).
 \label{eq:FC-rescaled-equations}
\end{equation}
Since $\widetilde{U}(0)=0$ and $D\widetilde{U}(0)=0$, it follows from \eqref{eq:FC13} that
\begin{equation}
 \|\widetilde{U}\|_{C^{1,1}(B_4)}\le C.
 \label{eq:FC-reference-bound}
\end{equation}
For $w=U_*-\widetilde{U}$, \eqref{eq:FC7} and \eqref{eq:FC-distance-scaling}
give
\begin{equation}
 \|w\|_{L^\infty(B_2)}
 \le Cr^{\alpha-2\theta}(1+|\log r|)\le Cd^{\alpha/2}.
 \label{eq:FC15}
\end{equation}
Here the logarithm is absorbed because
$(\alpha-2\theta)/(1+\theta)>\alpha/2$ by \eqref{eq:FC-theta-small}.
Using the common $C^{1,\beta}$ bound in \eqref{eq:FC3}, we also have
\begin{equation}
 [Dw]_{C^\beta(B_2)}
 =s^{\beta-1}
       [D(u_f-\widetilde u)]_{C^\beta(B_{2s}(x))}
 \le Cd^{\beta-1}.
 \label{eq:FC-gradient-seminorm}
\end{equation}
Set
\begin{equation}
 \gamma=1-\frac{\alpha}{4n}.
 \label{eq:FC-interpolation-exponents}
\end{equation}
Then $\gamma>1-\alpha/[n(2+\alpha)-1]$. Since
$1-\beta=\alpha^2/[64n^2(n+1)]<\alpha^2/(32n)$, we have
\[
 \beta-\gamma\ge\frac\alpha{8n},\qquad
 (\alpha/2)(\beta-\gamma)-(1-\beta)(1+\gamma)
       \ge\frac{\alpha^2}{16n}-2(1-\beta)>0.
\]
By \eqref{eq:FC15}--\eqref{eq:FC-gradient-seminorm} and the
interpolation inequality \cite[Proposition~2.2]{LLZ},
\begin{equation}
 \begin{aligned}
 \|w\|_{C^{1,\gamma}(B_1)}
 &\le C\|w\|_{L^\infty(B_2)}\\
 &\quad+C\|w\|_{L^\infty(B_2)}^{(\beta-\gamma)/(1+\beta)}
             [Dw]_{C^\beta(B_2)}^{(1+\gamma)/(1+\beta)}\le C.
 \end{aligned}
 \label{eq:FC-interpolation}
\end{equation}
Equations~\eqref{eq:FC-interpolation} and
\eqref{eq:FC-reference-bound} give a uniform $C^{1,\gamma}(B_1)$
bound for $U_*$. The rescaled right-hand side satisfies
\begin{equation}
 \lambda\le f(x+s\,\cdot)\le\Lambda;\qquad
 [f(x+s\,\cdot)]_{C^\alpha(B_2)}\le F_0s^\alpha.
 \label{eq:FC15b}
\end{equation}
The interior estimate \cite[Theorem~1.2]{LLZ}, applied on $B_1$,
therefore gives
\[
 |D^2U_*(0)|\le C,\qquad
 D^2u_f(x)=D^2U_*(0).
\]
The constants are uniform in $p,r,x$, so \eqref{eq:FC5} holds in
the fixed collar $d(x)\le(c_*/4)(r_0/2)^{1+\theta}$.
On the remaining compact set, the right-hand side of
\eqref{eq:PIA11} is uniformly bounded. This proves the proposition.
\end{proof}

\section{Proofs of the main theorems}\label{sec:main-proof}

In this section, we use the uniform estimates on smooth approximations to
construct the $C^{2,\alpha}$ solution and identify it with the
continuous Bedford--Taylor solution.

\begin{proof}[Proof of Theorems~\ref{thm:main} and~\ref{thm:existence}]
Let $(\Omega_j,\Phi_j,f_j)$ be the approximations of
Lemma~\ref{lem:smooth-data}. By \cite{CKNS}, there are smooth strictly
plurisubharmonic solutions $u_j$ on $\Omega_j$ with right-hand sides
$f_j$ and boundary values $\Phi_j|_{\partial\Omega_j}$.

On each $\Omega_j$, Proposition~\ref{prop:gradient-bootstrap}
gives \eqref{eq:FC3} for every smooth density satisfying \eqref{eq:FC2}.
Propositions~\ref{prop:finite-closure} and~\ref{prop:elliptic-finish}
then give $\|u_j\|_{C^{2,\alpha}(\bar\Omega_j)}\le C$.
The constant is independent of $j$ by the uniform estimates in
Subsection~\ref{sec:analytic-uniformity} and the common $C^3$ bounds
and third-derivative moduli used in the boundary estimates.

By Lemma~\ref{lem:endpoint-limit}, a subsequence of the
extensions converges to $v\in C^{2,\alpha}(\bar\Omega)$ satisfying
\eqref{intro:dirichlet} and \eqref{eq:SE2}.
Since $\det(v_{i\bar j})=f\ge\lambda>0$, its complex Hessian is
positive definite up to the boundary. If $u$ is a continuous
Bedford--Taylor solution with the same data, the comparison principle
\cite[Theorem~A]{BT} gives $u=v$.
\end{proof}

On the fixed domain of Theorem~\ref{thm:main}, the map
\[
 C^{2,\alpha}(\bar\Omega)\ni v
 \longmapsto\bigl(\det(v_{i\bar j}),\,v|_{\partial\Omega}\bigr)
\]
is continuously differentiable on the open set where
$(v_{i\bar j})>0$ on $\bar\Omega$.
At a solution $u$ of \eqref{intro:dirichlet}, its derivative
is $h\mapsto(fL_uh,h|_{\partial\Omega})$. This is an isomorphism by the linear Dirichlet theory and the global
Schauder estimate \cite[Theorems~6.6 and~6.14]{GT}.
By the implicit function theorem \cite[Theorem~17.6]{GT}
and uniqueness, the solution map is locally $C^1$ and hence locally
Lipschitz from
$C^\alpha(\bar\Omega)\times C^{2,\alpha}(\partial\Omega)$ into
$C^{2,\alpha}(\bar\Omega)$.

\begin{proof}[Proof of Corollary~\ref{cor:higher-regularity}]
Since $k\ge1$, Theorem~\ref{thm:main} with $\alpha=\sigma$
gives a $C^{2,\sigma}$ bound. The determinant lower bound gives
uniform ellipticity. Applying \cite[Proposition~7.1.4]{Han}
with derivative order $k$ and exponent $\sigma$, we obtain
\[
 u\in C^{k+2,\sigma}(\bar\Omega).
\]%

The norm bounds follow from the Schauder estimates
\cite[Theorem~6.6 and Section~6.4]{GT}. At step $j=1,\ldots,k$, we differentiate
tangentially. The coefficients and right-hand sides are bounded in
$C^{j-1,\sigma}$ by the $C^{j+1,\sigma}$ bound for $u$, the
$C^{j+2,\sigma}$ bounds for the boundary charts and values, and the
$C^{j,\sigma}$ bound for $f$. These estimates control derivatives involving at least one
tangential direction. The pure normal derivative is recovered from
the equation and uniform ellipticity, as in \cite[Theorem~6.19]{GT}.
\end{proof}

\end{document}